\documentclass[preprint,12pt]{elsarticle}

\usepackage{amsmath}
\usepackage{amsfonts}
\usepackage{mathrsfs}
\usepackage{amssymb, color}
\usepackage[linkcolor=black,anchorcolor=black,citecolor=black]{hyperref}
\usepackage{graphicx}
\usepackage[body={16.45cm,21.9cm}]{geometry}
\usepackage{amssymb}

\providecommand{\U}[1]{\protect\rule{.1in}{.1in}}

\providecommand{\U}[1]{\protect \rule{.1in}{.1in}}
\newtheorem{theorem}{Theorem}[section]

\newtheorem{definition}[theorem]{Definition}

\newtheorem{lemma}[theorem]{Lemma}

\newtheorem{proposition}[theorem]{Proposition}
\newtheorem{remark}[theorem]{Remark}

\newenvironment{proof}[1][Proof]{\noindent \textbf{#1.} }{\  \rule{0.5em}{0.5em}}

\usepackage{lineno,hyperref}

\def \R{\mathbb{R}}

\def \hE {\hat{\mathbb E}}

 \def\cliprr{C_{b,Lip}(\mathbb{R})}

\def\cliprm{C_{b,Lip}(\mathbb{R}^m)}
\def\cliprm+1{C_{b,Lip}(\mathbb{R}^{m+1})}

\numberwithin{equation}{section}

\begin{document}

\begin{frontmatter}



\title{The Cox-Ingersoll-Ross process under volatility uncertainty}

\author{Bahar Akhtari\fnref{a}}
\address[a]{Department of Actuarial Science, Faculty of Mathematical Sciences, Shahid Beheshti University, Tehran, Iran}
\ead{b\_akhtari@sbu.ac.ir}


\author{Hanwu Li\corref{cor1}\fnref{b,c}}
\address[b]{Research Center for Mathematics and Interdisciplinary Sciences, Shandong University,Binhai Rd. 72, Qingdao, 266237, Shandong, China}
\ead{lihanwu@sdu.edu.cn}
\cortext[cor1]{Corresponding Author}

\address[c]{Frontiers Science Center for Nonlinear Expectations (Ministry of Education), Shandong University,
Binhai Rd. 72, Qingdao, 266237, Shandong, China}





\begin{abstract}
Due to the importance of the Cox-Ingersoll-Ross process in different areas of finance, a broad spectrum of studies and investigations on this model have been carried out. In case of ambiguity, we characterize it by applying the $G$-expectation theory and the associated $G$-Brownian motion. In this paper, we provide the existence and uniqueness of the solution of the Cox-Ingersoll-Ross process in the presence of volatility uncertainty. In addition, some properties of the solution are indicated, such as the regularity and strong Markov property. Besides, we calculate some moments of the CIR process using a generalization of the nonlinear Feynman-Kac theorem.
\end{abstract}




\end{frontmatter}

\section{Introduction}
Knight was one of the key economists of the 20th century who
distinguished "risk" and "uncertainty" in his book \cite{KFH}. He illustrated that the notion of Risk refers to situations under which every occurrence unambiguously is attributed by a definite probability distribution. In contrast, the notion of ``Knightian uncertainty" or "Ambiguity", due to imperfect knowledge, refers to the circumstances under which the best decision can not be made by well-defined quantitative probabilities. In light of the uncertainty and lack of objective probability measures, a set of priors is considered to guide decision-makers under which ambiguity can be mitigated. One  of the frequently used tools to study financial and economic problems facing ambiguity, especially volatility uncertainty, is  the $G$-expectation theory established by Peng \cite{P}. In this framework, a new kind of Brownian motion called the $G$-Brownian motion, as well as the associated $G$-It\^{o}'s calculus, has been developed. Roughly speaking, the $G$-expectation is a nonlinear expectation induced by a non-dominated family of probability measures, and the $G$-Brownian motion is a continuous process starting from the origin with stationary and independent increments. It is worth pointing out that different from the classical case, the quadratic variation of the $G$-Brownian motion is not deterministic. 

The $G$-expectation has attracted plenty of attention both in theoretical analysis and applications. Gao \cite{G} considered the stochastic differential equations driven by a $G$-Brownian motion with Lipschitz coefficients ($G$-SDEs for short) and the strong Markov property for $G$-SDEs has been studied in \cite{HJL}. Then, Bai and Lin \cite{BL} get the existence and uniqueness result for $G$-SDEs with integral Lipschitz coefficients. Epstein and Ji \cite{EJ13} formulated a model of utility for a continuous time framework that captures the decision maker's concern with both drift and volatility uncertainty using $G$-expectation theory and in \cite{EJ14}, they derived arbitrage-free pricing rules and studied equilibrium in a representative agent economy with sequential security markets. It is worth pointing out that without using G-expectation theory, a version of the fundamental theorem of asset pricing in a continuous-time financial market under model uncertainty is proved in Biagini et al. \cite{SBKN} by assumption that each $P \in \mathcal{P}$ admits a martingale measure. Then optimal superhedging strategies for general contingent claims are provided and then by establishing the related duality, the superhedging price in terms of martingale measures
is derived. 


A diffusion process driven by $G$-Brownian is considered in H\"olzermann \cite{H3} to model the instantaneous forward rate. Then a drift condition consisting of several equations and several market prices is obtained for constructing an arbitrage-free term structure model which is robust with respect to the volatility. Also, it is shown that affine structure is preserved through some examples such as Ho-Lee term structure. In \cite{H1}, the Hull-White model under volatility uncertainty is studied under which it is shown that the classical martingale modeling does not work and the structure model should be updated by an adjustment factor to find an arbitrage-free term structure. Also, robust pricing formula for interest rate derivatives in fixed income markets under volatility uncertainty is discussed in \cite{H2}.

Akhtari et al. \cite{ABMO} by considering a general term structure of interest rates under model uncertainty, present a novel theorem to obtain the conditional $G$-expectation of a discounted payoff as the limit of $C^{1,2}$ solutions of some regularized PDEs. This result provides the minimal initial capital to discounted superhedging price with no-arbitrage assumption, too.

The Cox-Ingersoll-Ross model first was introduced by Cox et al. \cite{CIRS1, CIRS} in 1985 to formulate an equilibrium simple one-factor model for the evolution of interest rates. The model includes a square-root term in the diffusion coefficient of the instantaneous short-rate movements and does not leave the set of non-negative real numbers and also due to the analytical tractability, possesses a significant role in the finance/insurance area. The empirical implications of the CIR process of the term structure of interest rates advanced by Brown et al. \cite{BD} which show that the model can be fitted quite well to historical data and used in the different securities markets.  Furthermore, Brigo et al. \cite{BA} considered the two-dimensional shifted square-root diffusion for interest-rate and credit derivatives with positive stochastic intensity and then by proposing a new implicit Euler scheme while preserving the positivity and using Monte Carlo simulation, acquired the valuation of credit default swap (CDS for short). 
The forecasting of future expected interest rates, through Vasicek and CIR models, from observed empirical data is achieved in \cite{OMB} whereas the analytical tractability of the process is preserved.
The readers can refer to \cite{BM, MY} to see more discussions about the CIR process in the world of certainty. 
Fadina et al. \cite{FNS} developed a nonlinear affine process under parameter uncertainty which leads to a tractable model for Knightian uncertainty and discussed nonlinear versions of Vasicek and  Cox-Ingersoll-Ross (CIR for short) process which are suitable for modeling interest rates under ambiguity.




To the best of our knowledge, no systematic study concerning the CIR process under $G$-Brownian motion has been made. In the present paper, we aim to extend the classical CIR process to the case of ambiguity which is more complex than the classical case. Mathematically, we study the following type of $G$-SDEs
\begin{align}\label{mainGSD}
\begin{cases}
dX_t = (2\beta_1 X_t+\delta^1_t)dt + (2\beta_2 X_t+\delta^2_t)d\langle B\rangle_t+ \sigma(X_t)dB_t, \quad t \geq 0, \\[2mm]
X_0 = x_0,
\end{cases}
\end{align}
where $x_0$, $-\beta_1$, $-\beta_2$ are positive constants, $\delta^i$, $i=1,2$ are processes satisfying some appropriate integrability conditions and $\sigma$ defined on $[0,\infty)$ is a continuous function satisfying the $\frac{1}{2}$-H\"{o}lder condition 
\begin{displaymath}
|\sigma(x)-\sigma(y)|\leq K\sqrt{|x-y|}, \quad x,y \in [0,\infty).
\end{displaymath}
Motivated by \cite{DD}, we construct the solution by an approximation method. The main difficulty is to prove the convergence in an appropriate norm, which can be solved using some technical results in \cite{GR}. Although the approximate sequence may be negative, the limiting process is non-negative quasi-surely. 

The paper is organized as follows. In Section 2, we give some preliminaries in sublinear expectation theory. In Section 3, the existence and uniqueness of the non-negative solution for the $G$-SDE \eqref{mainGSD} are investigated. Section 4 is devoted to the regularity of the solution with respect to the starting time and original value.  Based on this result and by assumption that $\delta^i$, $i=1,2$ are deterministic functions, the strong Markov property is established. Section 5 provides some moment estimates for the solution $X$ which can be obtained by using a generalization of the nonlinear Feynman-Kac theorem for $\frac{1}{2}$-H\"{o}lder diffusion coefficient.

\section{Mathematical Preliminaries}
In this section, we review some of the essential tools in sublinear expectation theory. The reader can see more details in \cite{DHP11,HWZ, P0, P}.

Let $\Omega_T=C_{0}([0,T];\mathbb{R})$, the space of
real-valued continuous functions with $\omega_0=0$, be endowed
with the supremum norm and 
let  $B$ be the canonical process. Set
\[
L_{ip} (\Omega_T):=\{ \varphi(B_{t_{1}},...,B_{t_{n}}):  \ n\in\mathbb {N}, \ t_{1}
,\cdots, t_{n}\in\lbrack0,T], \ \varphi\in C_{b,Lip}(\mathbb{R}^{ n})\},
\]
where $C_{b,Lip}(\mathbb{R}^{ n})$ denotes the set of bounded Lipschitz functions on $\mathbb{R}^{n}$.
 We fix a sublinear and monotone function $G:\mathbb{R}\rightarrow\mathbb{R}$ defined by
\begin{displaymath}
G(a):=\frac{1}{2}(\overline{\sigma}^2a^+-\underline{\sigma}^2a^-), \quad \textup{for} \quad a \in \R,
\end{displaymath}
where $0\leq \underline{\sigma}^2<\overline{\sigma}^2$. The related $G$-expectation on $(\Omega_T, L_{ip}(\Omega_T))$ can be constructed in the following way. Assume that $\xi\in L_{ip}(\Omega_T)$ can be represented as
    \begin{displaymath}
    	\xi=\varphi(B_{{t_1}}, B_{t_2},\cdots,B_{t_n}).
\end{displaymath}
Then, for $t\in[t_{k-1},t_k)$, $k=1,\cdots,n$, set
\begin{displaymath}
	\hat{\mathbb{E}}_{t}[\varphi(B_{{t_1}}, B_{t_2},\cdots,B_{t_n})]=u_k(t, B_t;B_{t_1},\cdots,B_{t_{k-1}}),
\end{displaymath}
where $u_k(t,x;x_1,\cdots,x_{k-1})$ is a function of $(t,x)$ parameterized by $(x_1,\cdots,x_{k-1})$ such that it solves the following fully nonlinear PDE defined on $[t_{k-1},t_k)\times\mathbb{R}$:
\begin{displaymath}
	\partial_t u_k+G(\partial_x^2 u_k)=0
\end{displaymath}
with terminal conditions
\begin{displaymath}
	u_k(t_k,x;x_1,\cdots,x_{k-1})=u_{k+1}(t_k,x;x_1,\cdots,x_{k-1},x), \ k<n
\end{displaymath}
and $u_n(t_n,x;x_1,\cdots,x_{n-1})=\varphi(x_1,\cdots,x_{n-1},x)$. Hence, the $G$-expectation of $\xi$ is $\hat{\mathbb{E}}_0[\xi]$ and for simplicity, we always omit the subscript $0$. The triple $(\Omega_T, L_{ip}(\Omega_T),\hat{\mathbb{E}})$ is called the $G$-expectation space.


Define $\Vert\xi\Vert_{L_{G}^{p}}:=(\hat{\mathbb{E}}[|\xi|^{p}])^{1/p}$ for $\xi\in L_{ip}(\Omega_T)$ and $p\geq1$.   The completion of $L_{ip} (\Omega_T)$ under this norm  is denote by $L_{G}^{p}(\Omega_T)$. For all $t\in[0,T]$, $\hat{\mathbb{E}}_t[\cdot]$ is a continuous mapping on $L_{ip}(\Omega_T)$ w.r.t the norm $\|\cdot\|_{L_G^1}$. Hence, the conditional $G$-expectation $\mathbb{\hat{E}}_{t}[\cdot]$ can be
extended continuously to the completion $L_{G}^{1}(\Omega_T)$. Denis, Hu and Peng \cite{DHP11} prove that the $G$-expectation has the following representation.

\begin{theorem}[\cite{DHP11}]
\label{the1.1}  Let $\mathcal{B}(\Omega_T)$ be the Borel $\sigma$-algebra of $\Omega_T$. Then, there exists a weakly compact set
$\mathcal{P}$ of probability
measures on $(\Omega_T,\mathcal{B}(\Omega_T))$, such that
\[
\hat{\mathbb{E}}[\xi]=\sup_{P\in\mathcal{P}}\mathrm{E}_{P}[\xi] \text{ for all } \xi\in  {L}_{G}^{1}{(\Omega_T)}.
	\]
Moreover, $\mathcal{P}$ is called a set that represents $\hat{\mathbb{E}}$.
\end{theorem}

Let $\mathcal{P}$ be a weakly compact set that represents $\hat{\mathbb{E}}$.
For this $\mathcal{P}$, we define the capacity%
\[
c(A):=\sup_{P\in\mathcal{P}}P(A),\ A\in\mathcal{B}(\Omega_T).
\]
A set $A\in\mathcal{B}(\Omega_T)$ is called polar if $c(A)=0$.  A
property holds $``quasi$-$surely"$ $(q.s.)$ if it holds outside a
polar set. In the following, we do not distinguish two random variables $X$ and $Y$ if $X=Y$, $q.s.$

\begin{definition}\label{def2.6} 
Let $M_{G}^{0}(0,T)$ be the collection of processes in the
following form: for a given partition $\pi
_{T}:=\{t_{0},\cdot\cdot\cdot,t_{N}\}$ of $[0,T]$,
\begin{equation}\label{123}
\eta_{t}(\omega)=\sum_{j=0}^{N-1}\xi_{j}(\omega)\mathbf{1}_{[t_{j},t_{j+1})}(t),
\end{equation}
where $\xi_{i}\in L_{ip}(\Omega_{t_{i}})$, $i=0,1,2,\cdot\cdot\cdot,N-1$. For each
$p\geq1$ and $\eta\in M_G^0(0,T)$ let $\|\eta\|_{H_G^p}:=\{\hat{\mathbb{E}}[(\int_0^T|\eta_s|^2ds)^{p/2}]\}^{1/p}$, $\Vert\eta\Vert_{M_{G}^{p}}:=(\mathbb{\hat{E}}[\int_{0}^{T}|\eta_{s}|^{p}ds])^{1/p}$ and denote by $H_G^p(0,T)$,  $M_{G}^{p}(0,T)$ the completion
of $M_{G}^{0}(0,T)$ under the norm $\|\cdot\|_{H_G^p}$, $\|\cdot\|_{M_G^p}$ respectively.
\end{definition}

\begin{definition}
For each $\eta\in M_G^0(0,T)$ of the form \eqref{123}, we define the linear mappings $I,L:M_G^0(0,T)\rightarrow L_G^p(\Omega_T)$ as the following:
\begin{align*}
&I(\eta):=\int_0^T\eta_s dB_s=\sum_{j=0}^{N-1}\xi_j(B_{t_{j+1}}-B_{t_j}),\\
&L(\eta):=\int_0^T\eta_s d\langle B\rangle_s=\sum_{j=0}^{N-1}\xi_j(\langle B\rangle_{t_{j+1}}-\langle B\rangle_{t_j}).
\end{align*}
Then $I$ and $L$ can be continuously extended to $H_G^p(0,T)$ and $M_G^p(0,T)$, respectively.
\end{definition}

We have the following estimate for $G$-It\^{o}'s integral, which can be regarded as the B-D-G inequality under $G$-expectation.

\begin{proposition}[\cite{HJPS2}]\label{BDG}
If $\eta\in H_G^{\alpha}(0,T)$ with $\alpha\geq 1$ and $p\in(0,\alpha]$, then we get
\begin{displaymath}
\sup_{u\in[t,T]}|\int_t^u\eta_s dB_s|^p\in L_G^1(\Omega_T),
\end{displaymath}
and
\begin{displaymath}
\underline{\sigma}^p c_p\hat{\mathbb{E}}_t[(\int_t^T |\eta_s|^2ds)^{p/2}]\leq
\hat{\mathbb{E}}_t[\sup_{u\in[t,T]}|\int_t^u\eta_s dB_s|^p]\leq
\overline{\sigma}^p C_p\hat{\mathbb{E}}_t[(\int_t^T |\eta_s|^2ds)^{p/2}].
\end{displaymath}
\end{proposition}

Let $S_G^0(0,T)=\{h(t,B_{t_1\wedge t}, \ldots,B_{t_n\wedge t}):t_1,\ldots,t_n\in[0,T],h\in C_{b,Lip}(\mathbb{R}^{n+1})\}$. For $p\geq 1$ and $\eta\in S_G^0(0,T)$, set $\|\eta\|_{S_G^p}=\{\hat{\mathbb{E}}[\sup_{t\in[0,T]}|\eta_t|^p]\}^{1/p}$. Denote by $S_G^p(0,T)$ the completion of $S_G^0(0,T)$ under the norm $\|\cdot\|_{S_G^p}$. 

In the following, we give the characterization for the space $M_G^p(0,T)$, which is obtained in \cite{HWZ}. Set $\mathcal{F}_t = \mathcal{B}(\Omega_t)$ and $\mathbb{F} = (\mathcal{F}_t)_{t \in [0, T]}$, where $\Omega_t:=\{\omega_{\cdot\wedge t},\omega\in\Omega_T\}$. We define the following distance on $[0, T] \times \Omega_T$
\begin{align*}
\rho((t,\omega) ,(t',\omega')) = |t-t'| + \max_{0 \leq u \leq T}|\omega_u - \omega'_u|
\end{align*}
for $t,t'\in [0, T]$ and $\omega, \omega' \in \Omega_T$. Remember that a process $(\eta_t)_{t \in [0, T]}$ is said to be progressively measurable if its restriction
on $[0, t] \times \Omega_t$ is $\mathcal{B}([0, t]) \bigotimes \mathcal{F}_t$-measurable for every $t$. 

\noindent For any $p \geq 1$, set
\begin{align*}
\mathbb{M}^p(0, T) = \Big\{\eta\,\, \text{is progressively measurable on}\, [0, T] \times \Omega_T \,\,\text{and}\,\,  \hE\Big[\int_0^T |\eta_t|^p \,dt \Big] < \infty \Big\},
\end{align*}
and the corresponding capacity
\begin{align*}
\hat{c}(A) =\frac{1}{T} \hE\big[\int_0^T I_A(t, \omega)\,dt \big], \textup{for each progressively measurable set 
$A \subset [0, T] \times \Omega_T$.}
\end{align*}

\begin{definition}[\cite{HWZ}]
A progressively measurable process $\eta : [0, T ] \times \Omega_T \rightarrow  \R$ is called quasi-continuous $(q.c.)$, if for each $\varepsilon > 0$, there exists a progressively measurable open set $G$ in $[0, T ] \times \Omega_T$ such that $\hat{c}(G) < \varepsilon$ and $\eta|_{G^c}$ is continuous.
\end{definition}
\begin{definition}[\cite{HWZ}]
We say that a progressively measurable process
$\eta:[0, T] \times \Omega_T \rightarrow \R$ has a
quasi-continuous version if there exists a quasi-continuous process $\eta'$ such that $\hat{c}(\{\eta \ne \eta'\}) = 0$.
\end{definition}
\begin{theorem}[\cite{HWZ}]
For any $p \geq 1$, we have
\begin{align*}
M^p_G(0, T) = \Big\{\eta \in \mathbb{M}^p(0, T): &\lim_{N \rightarrow \infty} \hE\Big[\int_0^T |\eta_t|^p {\bf{1}}_{\{|\eta_t|>N\}}\,dt \Big] =0 \,
\text{and} \\ 
&\eta\,\, \text{has a} 
\text{quasi-continuous version}\Big\}.
\end{align*}
\end{theorem}

Let $L^0(\Omega_T)$ be the space of all $\mathcal{B}(\Omega_T)$-measurable real functions, $B_b(\Omega_T)$ all bounded functions in $L^0(\Omega_T)$ and $C_b(\Omega_T)$ all continuous functions in $B_b(\Omega_T)$. 
We set, for $p > 0$, the following spaces: 
\begin{align*}
\mathbb{L}^p:= \{X \in L^0(\Omega_T) : \hE[|X|^p]=\sup_{P \in \mathcal{P}} \mathrm{E}_P[|X|^p] < \infty \}
\end{align*}
and
\begin{align*}
\mathbb{L}_c^p:= \{X \in C_b(\Omega_T) : \hE[|X|^p]=\sup_{P \in \mathcal{P}} \mathrm{E}_P[|X|^p] < \infty \}.
\end{align*}
\begin{theorem}
Let $\mathcal{P}$ be weakly compact and let $\{X_n\}_{n=1}^{\infty} \subset \mathbb{L}_c^1$ and $X \in \mathbb{L}^1$ be such that $X_n \downarrow X$, q.s., as $n \rightarrow \infty$. Then $\hE[X_n] \downarrow \hE[X]$. 
\end{theorem}

\section{$G$-SDEs with non-Lipschitz diffusion coefficient} 
In this section, we study the existence and uniqueness result for $G$-SDE \eqref{mainGSD}, whose solution can be regarded as the CIR process under $G$-expectation. Before investigating \eqref{mainGSD}, let us first consider the uniqueness for a more general $G$-SDE of the following type:
\begin{equation}\label{GBessel}
dX_t = \mu(t,X_t)\,dt + b(t,X_t)\,d\langle B\rangle_t+\sigma(t,X_t)\,dB_t, \quad t \in [0, T],
\end{equation}
with $X_0 = x_0 \in \R$,  where $\mu,b,\sigma:[0,T]\times \Omega\times \R\rightarrow \R$. We suppose that the coefficients satisfy the conditions
\begin{equation}\begin{split}\label{con}
&|\mu(t,\omega,x) - \mu(t,\omega,y)| + |b(t,\omega,x)-b(t,\omega,y)|\leq K|x-y|,\\[2mm]
&|\sigma(t,\omega,x)-\sigma(t,\omega,y)|\leq h(|x-y|),
\end{split}\end{equation}
for any $t\in[0,T]$, $\omega\in \Omega$ and $x,y\in\mathbb{R}$, where $K$ is a positive constant and $h:[0,\infty)\rightarrow [0,\infty)$ is a strictly increasing function with $h(0)=0$ and
\begin{align}\label{hh}
\int_{(0,\varepsilon)}h^{-2}(x)dx=\infty
\end{align}
for any $\varepsilon>0$.
\begin{remark}\label{remm}
$X$ is called the solution of equation \eqref{GBessel} if $X\in M_G^1(0,T)$ and it satisfies this equation with $b(.,X_.) \in M^{1}_G(0,T)$ and $\sigma(.,X_.) \in M^{2}_G(0,T)$. 
\end{remark}

\begin{proposition}\label{prop1}
Let $X^i\in M_G^1(0,T)$ be the solution of equation \eqref{GBessel}. 
Then, we have $X_t^1 = X_t^2$ in $L^1_G(\Omega)$ sense for any $t \in [0, T]$.
\end{proposition}

\begin{proof}
By  the proof of Proposition 2.13 in Karatzas and Shreve \cite{KS}, there exists a strictly decreasing sequence $\{a_n\}_{n=0}^\infty$ with $a_0=1$, $\lim_{n\rightarrow\infty}a_n=0$ and $\int_{a_n}^{a_{a-1}}h^{-2}(x)dx=n$, for any $n\geq 1$. For any fixed $n\geq 1$, there exists a continuous function $\rho_n$ defined on $\mathbb{R}$ whose support is in $(a_n,a_{n-1})$ such that for any $x>0$, $0\leq \rho_n(x)\leq \frac{2}{nh^2(x)}$ and $\int_{a_n}^{a_{n-1}}\rho_n(x)ds=1$. We define
\begin{align}\label{psin}
\psi_n(x):=\int_0^{|x|}\int_0^y\rho_n(u)du dy, \qquad x\in\mathbb{R}.
\end{align}
Then, for each $n\geq 1$, $\psi_n$ is twice continuously differentiable with $|\psi'_n(x)|\leq 1$. Furthermore, the sequence $\{\psi_n\}_{n=1}^\infty$ is nondecreasing and $\lim_{n\rightarrow\infty}\psi_n(x)=|x|$ for $x\in\mathbb{R}$.

Now, set $\hat{X}_t=X_t^1-X_t^2$, $\hat{\mu}_t = \mu(t, X^1_t) -  \mu(t, X^2_t)$, $\hat{b}_t = b(t,X_t^1)-b(t,X_t^2)$ and $\hat{\sigma}_t=\sigma(t,X_t^1)-\sigma(t,X_t^2)$. Applying the $G$-It\^{o}'s formula implies that
\begin{displaymath}
\psi_n(\hat{X}_t)= \int_0^t \psi'_n(\hat{X}_s)\hat{\mu}_s ds+
\int_0^t \psi'_n(\hat{X}_s)\hat{b}_s d\langle B\rangle_s
+\int_0^t \frac{1}{2}\psi^{''}_n(\hat{X}_s)\hat{\sigma}_s^2d\langle B\rangle_s
+\int_0^t \psi'_n(\hat{X}_s)\hat{\sigma}_s dB_s.
\end{displaymath}
 By Remark \ref{remm}, $\psi'_n(\hat{X}_s)\hat{\sigma}_s \in M^2_G(0, T)$ and so
\begin{equation}\label{0}
\hat{\mathbb{E}}[\int_0^t \psi'_n(\hat{X}_s)\hat{\sigma}_s dB_s]=0.
\end{equation}
The properties of $\psi_n$ yield that
\begin{displaymath}
\hat{\mathbb{E}}[\int_0^t \psi^{''}_n(\hat{X}_s)\hat{\sigma}_s^2d\langle B\rangle_s]\leq \frac{2\overline{\sigma}^2 t}{n}.
\end{displaymath}
Note that $\overline{\sigma}^2$ is a constant defined by $\overline{\sigma}^2:=\hat{\mathbb{E}}[B^2_1]$ which differs from the function ${\sigma}$. Combining the above analysis, we get
\begin{align}\label{psiito}
\hat{\mathbb{E}}[\psi_n(\hat{X}_t)]&\leq 
\hE[\int_0^t \psi'_n(\hat{X}_s)\hat{\mu}_s ds]+
\hE[\int_0^t \psi'_n(\hat{X}_s)\hat{b}_s d\langle B\rangle_s]+\frac{\overline{\sigma}^2 t}{n}\\
&\leq K(1+\overline{\sigma}^2) \int_0^t \hE[|\hat{X}_s|]ds+ \frac{\overline{\sigma}^2 t}{n}. \nonumber
\end{align}
Letting $n\rightarrow \infty$, by the monotone convergence property for $G$-expectation, we have that $\hat{\mathbb{E}}[|\hat{X}_t|] \leq K (1 +\overline{\sigma}^2) \int_0^t \hat{\mathbb{E}}[|\hat{X}_s|]ds$ and applying the Gronwall inequality, we obtain $\hat{\mathbb{E}}[|\hat{X}_t|]=0$ for any $t\in[0,T]$. Finally by this fact that 
\begin{align*}
\hat{\mathbb{E}}[\int_0^T |\hat{X}_t|dt] \leq \int_0^T \hat{\mathbb{E}}[|\hat{X}_t|]dt=0,
\end{align*}
we find that $X^1=X^2$ in $M_G^1$ sense and the desired result follows.
\end{proof}

Now, inspired by application in finance, we consider the following type of one-dimensional $G$-SDE whose values are taken from $[0, \infty)$:
\begin{equation}\label{GSDE}
dX_t= (2\beta_1 X_t+\delta^1_t)dt + (2\beta_2 X_t+\delta^2_t)d\langle B\rangle_t+ \sigma(X_t)dB_t, \qquad t \in [0,T],
\end{equation}
with $X_0 = x_0 \geq 0$. The parameters $\beta_1, \beta_2$ are two negative constants and $\delta^1 , \delta^2 :[0,T] \times \Omega  \rightarrow \R$ are two stochastic processes which belong to $M_G^2(0,T)$. 
Moreover, $\sigma : [0, \infty) \rightarrow \R$ is a continuous function vanishing at origin and satisfying the H\"{o}lder condition 
\begin{equation}\label{sigma}
|\sigma(x)-\sigma(y)|\leq K\sqrt{|x-y|}, \quad x,y \in [0, \infty),
\end{equation}
where $K$ is a positive constant. 

 From \eqref{sigma}, we can write 
$ |\sigma(x)|\leq K\sqrt{|x|}$ and
 we can find a positive constant $\widetilde{K}$ only depending on $K$ such that
\begin{equation}\label{sigma2}
|\sigma(x)|^2\leq \widetilde{K}(1 + |x|^2).
\end{equation}
This confirms that $\sigma$ holds linear growth condition. In the following of this section, $C$ always represents a constant depending on $T,\overline{\sigma},\underline{\sigma},x_0,K,\beta_1,\beta_2,k$ with $k = \hE[\int_0^t |\delta^1_s|^2 ds] \vee \hE[\int_0^t |\delta^2_s|^2 ds]$, which may vary from line to line.

\begin{theorem}\label{main}
The $G$-SDE \eqref{GSDE} has a unique quasi-surely positive solution $X\in M_G^1(0,T)$. More precisely, we have $X\in S_G^2(0,T)$.
\end{theorem}

\begin{proof}
The uniqueness of the solution is obtained by Proposition \ref{prop1}.
The proof of existence is similar to the one in Deelstra and Delbaen \cite{DD}, which is divided into the following four steps.

{\color{black} \textbf{Step 1.}}
For each $n\geq 1$, we denote by $\pi^n_T$ the partition of $[0,T]$ as $\pi^n_T = \{t_0^n,t^n_1,\cdots, t^n_{N_n}\}$ with $0=t^n_0<t^n_1<\cdots<t^n_{N_n}=T$ such that $\lim _{n \rightarrow \infty} h_n = 0$, where
\begin{align}\label{step}
h_n:=\sup_{0\leq i\leq N_n-1}|t^n_{i+1}-t^n_i|.
\end{align}
For simplicity, we always omit the index $n$. We define Euler's "polygonal" approximation
$X_n(t)$ as follows:
\begin{align*}
X_n(t)&=X_n(t_k)+2\beta_1 X_n(t_k)(t-t_k) + \int_{t_k}^t\delta^1_u\,du
+ 2\beta_2 X_n(t_k)(\langle B\rangle_t-\langle B\rangle_{t_k})\\
&+\int_{t_k}^t\delta^2_ud\langle B\rangle_u
+{\sigma}(X_n(t_k))(B_t-B_{t_k}),
\end{align*}
where $t\in(t_k,t_{k+1}]$, $k=0,\cdots,N-1$, and $X_n(0)=x_0$. Since the
approximation $X_n$ might leave $[0, \infty)$, so we define a projection function as
\begin{displaymath}
\widetilde{\sigma}(x)=\begin{cases}
\sigma(x), & x \in [0, \infty);\\[3mm]
0, & x \not\in [0, \infty);
\end{cases}
\end{displaymath}
and then update the approximation as 
\begin{align}\label{app0}
X_n(t)&=X_n(t_k)+2\beta_1 X_n(t_k)(t-t_k)+\int_{t_k}^t\delta^1_u\,du
+ 2\beta_2 X_n(t_k)(\langle B\rangle_t-\langle B\rangle_{t_k})\\ \nonumber
&+\int_{t_k}^t\delta^2_ud\langle B\rangle_u
+\widetilde{\sigma}(X_n(t_k))(B_t-B_{t_k}).
\end{align}
It is easy to check that $\widetilde{\sigma}$ satisfies
\begin{align}\label{barsigma}
|\widetilde{\sigma}(x)-\widetilde{\sigma}(y)|\leq \widetilde{K} \sqrt{|x-y|},\, \quad x,y \in \R,
\end{align}
and
\begin{align}\label{barsigma1}
|\widetilde{\sigma}(x)|^2 \leq \widetilde{K} \big(1 + |x|^2),\, \quad x \in \R.
\end{align}
If we set $\eta_n(t)=t_k$ for $t_k\leq t<t_{k+1}$, the above equation can be written as
\begin{align}\label{app1}
X_n(t)-X_n(\eta_n(t))& = \int_{t_k}^t (2\beta_1 X_n(\eta_n(u))+\delta^1_u) du +
\int_{t_k}^t (2\beta_2 X_n(\eta_n(u))+\delta^2_u) d\langle B\rangle_u\\ \nonumber
&+\int_{t_k}^t \widetilde{\sigma}(X_n(\eta_n(u)))dB_u.
\end{align}

In order to analyze the upper bound of $X_n(t)-X_n(\eta_n(t))$, we first inspect $\hat{\mathbb{E}}[|X_n(\eta_n(t))|^2]$. To do this, from \eqref{app0} and by an induction method, we obtain
\begin{equation}\label{app2}
X_n(t)=x_0+ \int_0^t (2\beta_1 X_n(\eta_n(u))+\delta^1_u)\,du 
+ \int_0^t (2\beta_2 X_n(\eta_n(u))+\delta^2_u) d\langle B\rangle_u+\int_0^t \widetilde{\sigma}(X_n(\eta_n(u)))dB_u.
\end{equation}
So for $t\in[t_k,t_{k+1})$
\begin{align*}
\hE[|X_n(\eta_n(t))|^2]\leq & 4\,\{x_0^2 + \hE [(\int_0^{t_k} (2\beta_1 X_n(\eta_n(u))+\delta^1_u) du)^2]\\
&+ \hE [(\int_0^{t_k} (2\beta_2 X_n(\eta_n(u))+\delta^2_u) d\langle B\rangle_u)^2]+\hE[|\int_0^{t_k}\widetilde{\sigma}(X_n(\eta_n(u)))dB_u|^2]\}.
\end{align*}
Simple calculation implies that
\begin{align*}
\hE[|X_n(\eta_n(t))|^2]\leq & C\{1+\hE[\int_0^t |X_n(\eta_n(u))|^2du]+\hE[|\int_0^{t_k}|\widetilde{\sigma}(X_n(\eta_n(u)))|^2dB_u] \}\\
\leq &C\{1+\int_0^t \hE[|X_n(\eta_n(u))|^2]du\},
\end{align*}
where we have applied the relation \eqref{barsigma1}.  So, 
applying the Gr{o}nwall inequality yield that
\begin{equation}\label{e1}
\hat{\mathbb{E}}[|X_n(\eta_n(t))|^2]\leq C e^{Ct}.
\end{equation}
On the other hand, we can represent $X_n(\eta_n(.))$ as
\begin{align*}
X_n(\eta_n(t))=\sum_{k=0}^{N-1} X_n(t_k) I_{[t_k,t_{k+1})}(t) + X_n(T)I_{\{t=T\}}, \qquad t \in [0, T],
\end{align*}
where $X_n(t_k)\in L^2_G(\Omega_{t_k})$ and $X_n(T) \in L^2_G(\Omega_{T})$. So, $X_n(\eta_n(.)) \in M_G^2(0,T)$ and consequently,  $X_n(\cdot)\in S_G^2(0,T)$ by \eqref{app2}.

 Now, we set $U_n(t)=|X_n(t)-X_n(\eta_n(t))|$. In view of  \eqref{app1} and \eqref{e1}, we have
\begin{equation}\label{e2}
\hE[U^2_n(t)] \leq C h_n,
\end{equation}
where $h_n$ is defined by \eqref{step}. 

{\color{black} \textbf{Step 2.}} In the following, we prove that
\begin{displaymath}
\lim_{n,n'\rightarrow\infty}\hat{\mathbb{E}}[\sup_{t\in[0,T]}|X_n(t)-X_{n'}(t)|]=0.
\end{displaymath}
For this purpose, we will use the method introduced by Yamada \cite{Y}. We first find a sequence of numbers $1=a_0>a_1>\cdots>a_m>0$ such that $\int_{a_{i}}^{a_{i-1}}(Kx)^{-1} dx=i$, $i=1,\cdots,m$. It is easy to check that $\lim_{m\rightarrow\infty}a_m=0$.

 For each fixed $m \in \mathbb{N}$, we define $\varphi_m : \R \rightarrow \R$ as $\varphi_m(x)=\Phi_m(|x|)$, where $\Phi_m\in C^2([0,\infty); \mathbb{R})$ with $\Phi_m(0)=0$, such that $\Phi^{''}_m\in C([0,\infty); \mathbb{R})$ and
\begin{align}\label{fi"}
\Phi^{''}_m(x)=\begin{cases}
\textrm{between } 0 \textrm{ and } \frac{2}{mKx}, &a_m<x<a_{m-1},\\
0, &\textrm{otherwise}
\end{cases}
\end{align}
and $\int_{a_m}^{a_{m-1}}\Phi_m^{''}(x)dx=1$. It is easy to check that $|\Phi'_m(x)|\leq 1$, for $x\geq 0$, and $|x|-a_{m-1}\leq \varphi_m(x)$.

Set $X_{n,n'}(t) = X_n(t)-X_{n'}(t)$, $X_{n,n'}^\eta(t) = X_n(\eta_n(t))-X_{n'}(\eta_{n'}(t))$ and
$\sigma_{n,n'}^\eta(t) = \widetilde{\sigma}(X_n(\eta_n(t)))-\widetilde{\sigma}(X_{n'}(\eta_{n'}(t)))$. By G-It\^{o}'s formula, we have
\begin{align*}
\varphi_m(X_{n,n'}(t))&=\int_0^t 2\beta_1\varphi'_m(X_{n,n'}(u))X_{n,n'}^\eta(u)\,du
+\int_0^t \varphi'_m(X_{n,n'}(u))\sigma_{n,n'}^\eta(u) dB_u\\
&+\int_0^t 
\big(
2\beta_2\varphi'_m(X_{n,n'}(u))X_{n,n'}^\eta(u) +
\frac{1}{2}\varphi^{''}_m(X_{n,n'}(u))(\sigma_{n,n'}^\eta(u))^2 \big) d\langle B\rangle_u.
\end{align*}
Recalling equation \eqref{e1} and $|\varphi'_m|\leq 1$, the process $\varphi'_m(\hat{X}_\cdot)\sigma_{n,n'}^\eta(\cdot)$ belongs to $M^2_G(0,T)$. Hence,
\begin{displaymath}
\hat{\mathbb{E}}[\int_0^t \varphi'_m(X_{n,n'}(u))\sigma_{n,n'}^\eta(u)\,dB_u]=0.
\end{displaymath}
The first term can be written as
\begin{align*}
\hat{\mathbb{E}}[\int_0^t 2\beta_1\varphi'_m(X_{n,n'}(u))X_{n,n'}^\eta(u)\,du]  &=  \hat{\mathbb{E}}[\int_0^t 2\beta_1\varphi'_m(X_{n,n'}(u))\big(X_{n}(\eta_{n}(u)) - X_{n}(u)\\
& + X_{n'}(u)-X_{n'}(\eta_{n'}(u))
+ X_{n}(u) - X_{n'}(u) \big)\,du]
\end{align*}
Due to the fact that $\varphi'_m(x)\geq 0$ if $x\geq 0$ and $\varphi'_m(x)\leq 0$ if $x\leq 0$ and since $\beta_1 < 0$, we have $\beta_1 x\varphi'_m(x)\leq 0$. So, in view of \eqref{e2} and remember that $|\varphi'_m|\leq 1$, we have 
\begin{align}\label{rre1}
\hat{\mathbb{E}}[\int_0^t 2\beta_1\varphi'_m(X_{n,n'}(u))X_{n,n'}^\eta(u)\,du]
\leq  2|\beta_1|\int_0^t (\hat{\mathbb{E}}[|U_n(u)|]+\hat{\mathbb{E}}[|U_{n'}(u)|])du\leq  C(h^{\frac{1}{2}}_n + h^{\frac{1}{2}}_{n'}).
\end{align}
In a similar manner, we have
\begin{align}\label{re2}
\hat{\mathbb{E}}[\int_0^t 2\beta_1\varphi'_m(X_{n,n'}(u))X_{n,n'}^\eta(u)\,d\langle B\rangle_u]
&\leq  C (h^{\frac{1}{2}}_n + h^{\frac{1}{2}}_{n'}).
\end{align}

In the sequel, we obtain an upper bound for the last term. Noting that $\varphi^{''}_m(x)\leq \frac{2}{mK|x|}$. By setting $\|\varphi^{''}_m\|:=\sup_x|\varphi^{''}_m(x)|$, we have that 
\begin{align}\label{re3}
&\hE[\int_0^t \varphi^{''}_m(X_{n,n'}(u))(\sigma_{n,n'}^\eta(u))^2 d\langle B\rangle_u] \\
\leq & C
\hat{\mathbb{E}}[\int_0^t \varphi^{''}_m(X_{n,n'}(u))\big|X_{n}(\eta_{n}(u)) - X_{n}(u)
 + X_{n'}(u)-X_{n'}(\eta_{n'}(u))+ X_{n}(u) - X_{n'}(u) \big|\,du] \nonumber \\
\leq &C\|\varphi^{''}_m\|\int_0^t (\hE[|U_n(u)|]+\hat{\mathbb{E}}[|U_{n'}(u)|])\,du
+ C\hat{\mathbb{E}}[\int_0^t \frac{2}{m|X_{n,n'}(u)|K}|X_{n,n'}(u)|\,du]\nonumber \\
\leq & C\|\varphi^{''}_m\|(h^{\frac{1}{2}}_n + h^{\frac{1}{2}}_{n'}) + \frac{C}{m}. \nonumber 
\end{align}
So, by \eqref{rre1}, \eqref{re2}, \eqref{re3} and recalling that $|x|-a_{m-1}\leq \varphi_m(x)$, we yield that
\begin{align*}
\hat{\mathbb{E}}[|X_n(t)-X_{n'}(t)|] &\leq a_{m-1}+\hat{\mathbb{E}}[\varphi_m(X_n(t)-X_{n'}(t))]\\
& \leq a_{m-1} + 
C\big(2+\|\varphi^{''}_m\|\big)(h^{\frac{1}{2}}_n + h^{\frac{1}{2}}_{n'})
+ \frac{C}{m}.
\end{align*}
For any given $\varepsilon>0$, let $m$ be large enough such that $a_{m-1}<\frac{\varepsilon}{3}$ and $\frac{C}{m}<\frac{\varepsilon}{3}$. For this fixed $m$, noting that by \eqref{fi"} we have that $\|\varphi^{''}_m\| \leq \frac{2}{ma_mK}$, $\|\varphi^{''}_m\|$ is bounded. So, we can find $n_0$ large enough such that for any $n,n'\geq n_0$, 
\begin{align*}
C\big(2+\|\varphi^{''}_m\|\big)(h^{\frac{1}{2}}_n + h^{\frac{1}{2}}_{n'})<\frac{\varepsilon}{3}.
\end{align*}
All the above analysis implies that, for any $n,n'\geq n_0$ and for any $t\in[0,T]$,
\begin{displaymath}
\hat{\mathbb{E}}[|X_n(t)-X_{n'}(t)|]\leq a_{m-1}+\hat{\mathbb{E}}[\varphi_m(X_n(t)-X_{n'}(t))]<\varepsilon,
\end{displaymath}
which indicates that $\{X_n\}_{n\geq 1}$ is a Cauchy sequence in $M_G^1(0,T)$ and since $M_G^1(0,T)$ is a Banach space, there exists a process $X\in M_G^1(0,T)$ such that
\begin{align}\label{norm1}
\lim_{n\rightarrow\infty}\hat{\mathbb{E}}[\int_0^T |X_n(t)-X(t)|\,dt]=0.
\end{align}
By the result in Step 1, we have
\begin{align*}
\lim_{n\rightarrow\infty}\hat{\mathbb{E}}[\int_0^T |X_n(\eta_n(t))-X(t)|dt]=0.
\end{align*}

Recalling the dynamics of $X_n$ and $X_{n'}$ and by the B-D-G inequality in Proposition \ref{BDG}, it is easy to check that
\begin{align*}
\hat{\mathbb{E}}[\sup_{t\in[0,T]}|X_{n,n'}(t)|]\leq &C\{\hat{\mathbb{E}}[\int_0^T |X_{n,n'}^\eta(u)|\,du]
+\hat{\mathbb{E}}[\sup_{t\in[0,T]}|\int_0^t \sigma_{n,n'}^\eta(u)\, dB_u|]\}\\
\leq & C\{\hat{\mathbb{E}}[\int_0^T |X_{n,n'}^\eta(u)|\,du]+ \big(\hat{\mathbb{E}}[\int_0^T |\sigma_{n,n'}^\eta(u)|^2\,du]\big)^{1/2}\}\\
\leq & C\{\hat{\mathbb{E}}[\int_0^T |X_{n,n'}^\eta(u)|\,du]+  \big(\hat{\mathbb{E}}[\int_0^T |X_{n,n'}^\eta(u)|\, du]\big)^{1/2}\}.
\end{align*}
 We finally get that $\{X_n\}_{n\geq 1}$ is a Cauchy sequence in $S^1_G(0,T)$, i.e.,
\begin{displaymath}
\lim_{n,n'\rightarrow\infty}\hat{\mathbb{E}}[\sup_{t\in[0,T]}|X_n(t)-X_{n'}(t)|]=0.
\end{displaymath}
Therefore, the limit $X\in S_G^1(0,T)$ satisfies
\begin{align}\label{normsup}
\lim_{n\rightarrow\infty}\hat{\mathbb{E}}[\sup_{t\in[0,T]}|X_n(t)-X(t)|]=0.
\end{align}
Letting $n$ approach infinity in equation \eqref{app2}, we obtain that
\begin{equation}\label{e4}
X(t)=x_0+\int_0^t (2\beta_1 X(u)+\delta^1_u)\,du+\int_0^t (2\beta_2 X(u)+\delta^2_u)\,d\langle B\rangle_u+\int_0^t \widetilde{\sigma}(X(u))\,dB_u.
\end{equation}


{\color{black} \textbf{Step 3.} } In the following, we aim to show that
\begin{equation}\label{S2}
\lim_{n\rightarrow\infty}\hat{\mathbb{E}}[\sup_{t\in[0,T]}|X_n(t)-X(t)|^2]=0.
\end{equation}
 For this purpose, a similar analysis as Proposition 2.2 in Gy\"{o}ngy and R\'{a}sonyi \cite{GR} is presented. 
Let $\delta>1$ and $\varepsilon>0$. There exists a continuous nonnegative function $\psi_{\delta\varepsilon}(x)$, $x\in[0,\infty)$, whose support is in $[\varepsilon/\delta,\varepsilon]$, has integral $1$ and satisfies
\begin{displaymath}
\psi_{\delta\varepsilon}(x)\leq \frac{2}{x\ln \delta}.
\end{displaymath}
Set
\begin{align}\label{Phi}
\phi_{\delta\varepsilon}(x):=\int_0^{|x|}\int_0^y \psi_{\delta\varepsilon}(z)dz dy, \quad x \in \mathbb{R}.
\end{align}
It is easy to check that for any $x\in\mathbb{R}$, $|x|\leq \phi_{\delta\varepsilon}(x)+\varepsilon$, $|\phi'_{\delta\varepsilon}|\leq 1$, and
\begin{align*}
\phi^{''}_{\delta\varepsilon}(x)&=\psi_{\delta\varepsilon}(|x|)
\leq \frac{2}{|x|\ln \delta}I_{[\varepsilon/\delta,\varepsilon]}(|x|).
\end{align*}
First, set $Y_n(t) = X(t) - X_n(t)$. So $|Y_n(t)|\leq \varepsilon + \phi_{\delta\varepsilon}(Y_n(t))$.
After applying $G$-It\^{o}'s formula for $\phi_{\delta\varepsilon}(Y_n(t))$, we have
\begin{align*}
|Y_n(t)|&\leq \varepsilon +\phi_{\delta\varepsilon}(Y_n(t))\\
&=\varepsilon + \int_0^t l_{\delta\varepsilon n}(s)\,ds + \int_0^t I_{\delta\varepsilon n}(s)\,d\langle B\rangle_s+\frac{1}{2}\int_0^t J_{\delta\varepsilon n}(s)\,d\langle B\rangle_s + M_{\delta\varepsilon n}(t),
\end{align*}
where
\begin{align*}
&l_{\delta\varepsilon n}(s)=2\beta_1\phi'_{\delta\varepsilon}(Y_n(s))(X(s)-X_n(\eta_n(s))),\\[2mm]
&I_{\delta\varepsilon n}(s)=2\beta_2\phi'_{\delta\varepsilon}(Y_n(s))(X(s)-X_n(\eta_n(s))),\\[2mm]
&J_{\delta\varepsilon n}(s)=\phi^{''}_{\delta\varepsilon}(Y_n(s))(\widetilde{\sigma}(X(s))-\widetilde{\sigma}(X_n(\eta_n(s))))^2,\\[2mm]
&M_{\delta\varepsilon n}(t) = \int_0^t \phi^{'}_{\delta\varepsilon}(Y_n(s))(\widetilde{\sigma}(X(s))-\widetilde{\sigma}(X_n(\eta_n(s))))\,dB_s.
\end{align*}
Recalling that $X(\cdot), X_n(\eta_n(\cdot))\in M^2_G(0,T)$ and $|\phi^{'}_{\delta\varepsilon}|\leq 1$, we have $\phi^{'}_{\delta\varepsilon}(Y_n(\cdot))(\widetilde{\sigma}(X(\cdot))-\widetilde{\sigma}(X_n(\eta_n(\cdot))))\in M_G^2(0,T)$. Hence, $M_{\delta\varepsilon n}$ is a symmetric $G$-martingale. By the construction of $\phi_{\delta\varepsilon}$ and the assumption on $\widetilde{\sigma}$, it is easy to check that
\begin{align*}
l_{\delta\varepsilon n}(s)\leq 2|\beta_1|(|Y_n(s)|+ U_n(s)), \quad
I_{\delta\varepsilon n}(s)\leq 2|\beta_2|(|Y_n(s)|+ U_n(s)),
\end{align*}
and
\begin{align*}
J_{\delta\varepsilon n}(s)\leq \frac{2K^2}{\ln \delta}+\frac{2K^2 \delta}{\varepsilon\ln \delta}U_n(s).
\end{align*}
The above analysis implies that
\begin{align*}
|Y_n(t)|&\leq \varepsilon + 2|\beta_1| \int_0^t |Y_n(s)| ds
+2|\beta_1| \int_0^t U_n(s)\,ds +  \int_0^t \big(2|\beta_2||Y_n(s)| + \frac{K^2}{\ln \delta}\big)\,d \langle B \rangle_s \\
&+ \int_0^t \big(2|\beta_2| + \frac{K^2 \delta}{\varepsilon \ln\delta}\big) U_n(s)\,d \langle B \rangle_s
+M_{\delta\varepsilon n}(t).
\end{align*}
We now set $Z_n(t)=\sup_{s\in[0,t]}|Y_n(s)|$. It follows from equation above
\begin{align*}
\hE[Z^2_n(t)] \leq &C\,\big(\varepsilon^2 + \hE \big( \int_0^t [Z_n(s)]\,ds \big)^2
+ \int_0^t \hE[U^2_n(s)]\,ds + \frac{1}{(\ln\delta)^2} \\
& 
+ \int_0^t \big(1+ \frac{\delta^2}{\varepsilon^2 (\ln\delta)^2}\big) \hE[U^2_n(s)]\,ds+\hE[\sup_{s\in[0,t]}|M_{\delta\varepsilon n}(s)|^2]\big).
\end{align*}
Noting that $|\phi^{'}_{\delta\varepsilon}|\leq 1$ and by the B-D-G inequality in Proposition \ref{BDG}, we obtain that
\begin{align*}
\hE[\sup_{s\in[0,t]}|M_{\delta\varepsilon n}(s)|^2]\leq &  C \hat{\mathbb{E}}[\int_0^t (\widetilde{\sigma}(X(s))-\widetilde{\sigma}(X_n(\eta_n(s))))^2 ds]\\
\leq & C \big({\hE}[\int_0^t |Y_n(s)|ds]+ \hE[\int_0^t U_n(s)ds]\big).
\end{align*}
We set $\varepsilon = \frac{-1}{\ln h_n}$ and $\delta = h^{-\frac{1}{3}}_n$ when $n\geq 2$. In view of \eqref{e2}, all the above analysis implies that
\begin{align*}
\hE[Z^2_n(t)] &\leq {C}\,\big(\hE \big( \int_0^t [Z_n(s)]\,ds \big)^2 + \hE[\int_0^t |Y_n(s)|ds] + \frac{1}{(\ln h_n)^2}
+ h^{1/3}_n\big).
\end{align*}
 Applying Lemma 3.2 in Gy\"{o}ngy and R\'{a}sonyi \cite{GR} with linear expectation $\mathrm{E}$ replaced by $G$-expectation $\hat{\mathbb{E}}$, there exists a positive constant $\hat{C}$ independent of $n$, such that
\begin{align*}
\hat{\mathbb{E}}[Z^2_n(T)]\leq \hat{C} \big(\frac{1}{(\ln h_n)^2}+h_n^{1/3} + \int_0^T \hat{\mathbb{E}}[|Y_n(s)|]ds\big).
\end{align*}
In view of \eqref{norm1}, we obtain that
\begin{equation}\label{norm2}
\lim_{n\rightarrow\infty}\hat{\mathbb{E}}[\sup_{t\in[0,T]}|X(t)-X_n(t)|^2]=0.
\end{equation}
Hence, the claim holds true. Recalling that we have shown in Step 1 that $X_n(.) \in S_G^2(0,T)$, \eqref{S2} implies that $X\in S_G^2(0,T)$.

{\color{black} \textbf{Step 4.}} Now we prove that $X$ is the non-negative solution to \eqref{GSDE}. The proof of non-negativity is similar to the one in Deelstra and Delbaen \cite{DD} and the main difference is that the results hold for each $P\in\mathcal{P}$, so we omit it. Note that since $X\in M_G^1$, so it has a quasi-continuous modification and clearly  $\sigma(X_\cdot)$ has a quasi-continuous modification, too. Also
\begin{align}\label{MG1}
\lim_{N\rightarrow\infty}\hat{\mathbb{E}}[\int_0^T |X_t| I_{\{|X_t|\geq N\}}dt]=0,
\end{align}
and by \eqref{sigma}, we have that
\begin{align*}
 |\sigma(X_t)|^2\leq K^2|X_t|
\end{align*}
for any $t \in [0, T]$. If $|\sigma(x)|\geq N$ for $N$ large enough, we get
\begin{align*}
\{|\sigma(X_t)|\geq N\}\subset \{|X_t|\geq (\frac{N}{K})^2\},
\end{align*}
and so
\begin{equation}\label{00}
\lim_{N\rightarrow\infty}\hat{\mathbb{E}}[\int_0^T |\sigma(X_t)|^2 I_{\{|\sigma(X_t)|\geq N\}}dt]=0,
\end{equation}
which shows $\sigma(X_.) \in M^2_G(0, T)$ and it makes the integral $\int_0^. {\sigma}(X(u))\,dB_u$ well-defined.   
Now, letting $n$ approach infinity in equation \eqref{app2} and noting the definition of $\widetilde{\sigma}$, we have
\begin{equation*}
X(t) = x_0+\int_0^t (2\beta_1 X(u)+\delta^1_u)\,du+\int_0^t (2\beta_2 X(u)+\delta^2_u)\,d\langle B\rangle_u+\int_0^t {\sigma}(X(u))\,dB_u,
\end{equation*}
which shows that the equation \eqref{GSDE} has a non-negative solution.
\end{proof}

\begin{remark}
Note that approximation \eqref{app0} might leave $[0, \infty)$ but constructing an efficient method for the $G$-SDE \eqref{GSDE} is beyond the scope of this paper. We just utilized the construction approach to prove the existence of the solution. 
\end{remark}

\begin{remark}
We first proved that the sequence $(X_n)_{n \geq 0}$ converges to $X$ under the norm $\|\cdot\|_{M_G^1}$, see \eqref{norm1}, and then this fact was established under the norm $\|\cdot\|_{S_G^1}$, see \eqref{normsup}. Note that the norm $\|\cdot\|_{S_G^1}$ is finer than the norm $\|\cdot\|_{M_G^1}$. Indeed, we utilize the former result to achieve the latter one. 
\end{remark}

\begin{remark}
Note that the solution of \eqref{GSDE} belongs to $S_G^p(0,T)$ for $p \geq 2$ by assumption that $\delta^i\in M_G^p(0,T)$, for $i = 1,2$.
\end{remark}

\section{Strong Markov Property}
In this section, we intend to investigate the \emph{strong Markov property} of the process \eqref{GSDE}. For this purpose, first, we present some useful estimates.
\begin{lemma}\label{pbound}
Let $p \geq 2$ and $\delta^i\in M_G^p(0,T)$, $i=1,2$. The unique solution $(X^{t, x}_s)_{s \in [t, T]}$ of the $G$-SDE \eqref{GSDE} with value $x$ at initial time $t$ belong to $S^p_G(t, T)$ and there exists a constant $C$, which may vary from line to line, depending on $p, K, \overline{\sigma}$, $\beta_1$, $\beta_2, T,k$ with $k = \hE[\int_0^t |\delta^1_s|^p ds] \vee \hE[\int_0^t |\delta^2_s|^p ds]$ such that 
\begin{align}\label{re1}
\hE[\sup_{s\in[0,t]}|X^{x}_s|^p] \leq C(1+ |x|^p),
\end{align}
and for all $x, y \in \mathbb{R}$ and $t,t' \in [0, T]$, we have
\begin{align*}
\hE|X^{x}_t - X^{y}_{t'}|^p \leq C\big(|x - y|^{p} + |x-y| + (1+|x|^p)|t - t'|^{p/2}\big),
\end{align*}
and for $\gamma \in [0, T-t]$, we have that
\begin{align}\label{Xx}
\hE[\sup_{s\in[t,t+\gamma]}|X^{t,x}_s - x|^p] \leq C(1+ |x|^p)\gamma^{\frac{p}{2}}.
\end{align}
\end{lemma}

\begin{proof}
Using the fact $|a+b|^r \leq \max\{1, 2^{r-1}\} (|a|^r + |b|^r)$, $r >0$, and also by the H\"{o}lder and B-D-G inequality in Proposition \ref{BDG}
\begin{align*}
\hE[\sup_{s\in[0,t]}|X^x_s|^p] \leq & C \big\{|x|^p +   \int_0^t \hE[|X^x_u|^p]\,du +\hE[\int_0^t |\delta^1_s|^p ds]\\
&+\hE[\int_0^t |\delta^2_s|^p ds]+
\int_0^t (1 + \hE[|X^x_u|^p])\,du \big\}.
\end{align*}
Hence 
\begin{align*}
1 + \hE [\sup_{s\in[0,t]}|X^x_s|^p] \leq & C (1 + |x|^p) + C \big(\int_0^t (1 + \hE [\sup_{s\in[0,u]}|X^x_s|^p])\, du \big).
\end{align*}
 The desired result is followed by applying the Gronwall lemma: 
\begin{align}\label{c1c2}
1 + \hE [\sup_{s\in[0,t]}|X^x_s|^p] \leq & C e^{Ct} (1 + |x|^p).
\end{align}

Now we deal with the second inequality. To this end, first, set $Y^{x,y}_s = X^{x}_s - X^{y}_s$ with the following representation:
\begin{equation*}
Y^{x,y}_s = x - y + 2\beta_1 \int_{0}^{s} (X^{x}_u - X^{y}_u)\,du + 2\beta_2 \int_{0}^{s} (X^{x}_u - X^{y}_u)\,d\langle B\rangle_u + \int_0^s (\sigma(X^{x}_u) - \sigma(X^{y}_u))\,dB_u
\end{equation*}
for $s \in [0, T]$. We can write $|Y^{x,y}_s|\leq \varepsilon + \phi_{\delta\varepsilon}(Y^{x,y}_s)$ where $\phi_{\delta\varepsilon}$ is defined by \eqref{Phi} with $\delta>1$ and $\varepsilon>0$. After applying $G$-It\^{o}'s formula for $\phi_{\delta\varepsilon}(Y^{x,y}_s)$, we have
\begin{align*}
|Y^{x,y}_s|\leq &\varepsilon +\phi_{\delta\varepsilon}(Y^{x,y}_s)\\
=&\varepsilon+|x-y|+ \int_0^s l_{\delta\varepsilon n}(u)\,du + \int_0^s I_{\delta\varepsilon n}(u)\,d\langle B\rangle_u\\
&+\frac{1}{2}\int_0^s J_{\delta\varepsilon n}(u)\,d\langle B\rangle_u + M_{\delta\varepsilon n}(s),
\end{align*}
with
\begin{align*}
&l_{\delta\varepsilon n}(u)=2\beta_1\phi'_{\delta\varepsilon}(Y^{x,y}_u)(X^{x}_u - X^{y}_u),\\[2mm]
&I_{\delta\varepsilon n}(u)=2\beta_2\phi'_{\delta\varepsilon}(Y^{x,y}_u)(X^{x}_u - X^{y}_u),\\[2mm]
&J_{\delta\varepsilon n}(u)=\phi^{''}_{\delta\varepsilon}(Y^{x,y}_u)(\sigma(X^{x}_u) - \sigma(X^{y}_u))^2,\\[2mm]
& M_{\delta\varepsilon n}(s) = \int_0^s \phi^{'}_{\delta\varepsilon}(Y^{x,y}_u)(\sigma(X^{x}_u) - \sigma(X^{y}_u))\,dB_u.
\end{align*}
Recall that for any $x\in\mathbb{R}$, $|x|\leq \phi_{\delta\varepsilon}(x)+\varepsilon$, $|\phi'_{\delta\varepsilon}|\leq 1$, and
$\phi^{''}_{\delta\varepsilon}(x) \leq \frac{2}{|x|\ln \delta}I_{[\varepsilon/\delta,\delta]}(|x|)$. Note that $M_{\delta\varepsilon n}$ is a $G$-martingale. Henceforth
\begin{align*}
\hat{\mathbb{E}}[|Y^{x, y}_s|] \leq & {\varepsilon} + \big(|x - y| +  2(\beta_1 + \overline{\sigma}^{2} \beta_2) \int_0^s \hE[|Y^{x,y}_u|]\,du +  s\overline{\sigma}^2 \frac{K^2}{\ln \delta}
 \big),
\end{align*}
by applying Gronwall inequality and letting $\delta \rightarrow \infty$ and $\varepsilon \rightarrow 0$, we can find a positive constant $C$ independent of $x,y, s$ such that 
\begin{align}\label{pone}
\hat{\mathbb{E}}[|X^{x}_s - X^y_s|] \leq & C\,|x - y|.
\end{align}

Let $p \geq 2$. Then by Proposition \ref{BDG}, we have
\begin{align}\label{xyt'p}
\hE[|X^{x}_{t} - X^{y}_{t}|^p] &\leq C \big(|x - y|^p +  \int_{0}^{t} \hE[|X^{x}_u - X^{y}_u|^p]\,du 
+\int_0^t\hE[|X^{x}_u - X^{y}_u|^{\frac{p}{2}}]\,du\big).
\end{align}
Applying the H\"{o}lder inequality yields that
\begin{align}\label{phalf}
\hE[|X^{x}_u - X^{y}_u|^{\frac{p}{2}}] \leq C\big(\hE[|X^{x}_u - X^{y}_u|] + \hE[|X^{x}_u - X^{y}_u|^p]\big).
\end{align} 
By \eqref{pone}, \eqref{xyt'p} and \eqref{phalf}, we get
\begin{align*}
\hE[|X^{x}_{t} - X^{y}_{t}|^p] &\leq C\big(|x - y|^p
+ |x - y|+\int_{0}^{t} \hE[|X^{x}_u - X^{y}_u|^p]\,du\big).
\end{align*}
Using the Gronwall inequality implies that 
\begin{align}\label{xy}
\hE[|X^{x}_{t} - X^{y}_{t}|^p] &\leq C \big(|x - y|^p + |x - y|\big).
\end{align}
Now for $t', t \in [0, T]$, $t' > t$, we can write
\begin{align}\label{txt'}
X^{x}_t - X^{x}_{t'} = &2\beta_1 \int_{t}^{t'} X^{x}_u \,du + 2\beta_2 \int_{t}^{t'} X^{x}_u\,d\langle B\rangle_u + \int_t^{t'} \sigma(X^{x}_u)\,dB_u\\
& + \int_{t}^{t'} \delta_u^1\,du + \int_{t}^{t'} \delta_u^2\,d\langle B\rangle_u. \nonumber
\end{align}
By H\"{o}lder inequality and Proposition \ref{BDG}, we obtain that
\begin{equation}\begin{split}\label{txt'}
&\hE[|X^{x}_t - X^{x}_{t'}|^p]\\ \leq & C \big( (t'-t)^{p-1} \int_{t}^{t'} \hE[|X^{x}_u|^p] \,du+
\hE[\big( \int_t^{t'} |\sigma(X^{x}_u)|^2\,du \big)^{\frac{p}{2}}]\\
& + (t' - t)^{p-1}(\hE[\int_{t}^{t'}|\delta_u^1|^pdu] +\hE[\int_t^{t'}|\delta_u^2|^p\,du])\big) \\
\leq & C\big( (t'-t)^{p-1} \int_{t}^{t'} \hE[|X^{x}_u|^p] \,du + (t' -t )^{p-1}
+ (t' - t)^{\frac{p}{2}-1} \int_t^{t'} (1 + \hE[|X^{x}_u|^p])\,du \big).
\end{split}\end{equation}
Plugging Equation \eqref{c1c2} into the above inequality indicates that
 \begin{align}\label{txt'}
 \hE[|X^{x}_t &- X^{x}_{t'}|^p]  \leq C(1+|x|^p)(t'-t)^{\frac{p}{2}}.
 \end{align}
By combining \eqref{xy} and \eqref{txt'}, the desired result follows.

 For the last result, simple calculation yields that 
\begin{align*}
&\hE[\sup_{s\in[t, t+\gamma]}|X^{t,x}_s - x|^p] \\
\leq  &C\big(\hE[
(\int_t^{t+\gamma} |X^{t, x}_s|ds)^p]+(\int_t^{t+\gamma}
\hE[|\sigma(X^{t, x}_s)|^2]\,ds )^{\frac{p}{2}-1}\big) \\
&+\hE[(\int_t^{t+\gamma}|\delta^1_s|ds)^p]+\hE[(\int_{t}^{t+\gamma}|\delta^2_s|ds)^p]\\
\leq & C\big\{\gamma^{p-1} (\int_t^{t+\gamma} \hE[|X^x_s|^p]\,ds +\int_t^{t+\gamma} \hE[|\delta^1_s|^p]\,ds+\int_t^{t+\gamma} \hE[|\delta^2_s|^p]\,ds)\\
&+ \gamma^{\frac{p}{2}-1} \int_t^{t+\gamma} (1 + \hE[|X^x_s|^p])\,ds \big\}.
\end{align*}
By \eqref{re1}, we finally obtain that 
\begin{align*}
\hE[\sup_{s\in[t, t+\gamma]}|X^{t,x}_s - x|^p] \leq C(1+|x|
^p)\gamma^{\frac{p}{2}}.
\end{align*}
\end{proof}

In the sequel, we consider a more general version of the $G$-SDE \eqref{GSDE} as
\begin{equation}\label{GSDE1}
dX^{t, \xi}_s= (2\beta_1 X^{t, \xi}_s+\delta^1_s)\,ds + (2\beta_2 X^{t, \xi}_s+\delta^2_s)\,d\langle B\rangle_s+ \sigma(X^{t, \xi}_s)\,dB_s, \quad s \geq t,
\end{equation}
 where $X^{t, \xi}_t = \xi$ with $\xi \in L^p_G(\Omega_t)$, for some $p \geq 1$, and $\xi \geq 0$ quasi-surely.

 \begin{remark}
Modifying the proof of Theorem \ref{main}, for any given non-negative $\xi\in L_G^p(\Omega_t)$ with some $p\geq 1$, the $G$-SDE \eqref{GSDE1} admits a non-negative unique solution $X^{t,\xi}\in S_G^p(t,T)$. Besides, similar regularity property as in Lemma \ref{pbound} still holds. More precisely, let $\xi, \xi' \in L^p_G(\Omega_t)$ be non-negative, $X^{t,\xi}_s$ and $X^{t,\xi'}_s$, $s \geq t$, be the solutions of the $G$-SDE \eqref{GSDE1} with $X^{t,\xi}_t = \xi$ and $X^{t,\xi'}_{t} = \xi'$, respectively. Then, we have
\begin{align}\label{randp}
\hE_t[\sup_{s\in[t,T]}|X^{t,\xi}_s|^p] \leq C(1+ |\xi|^p),
\end{align}
and for $s,s' \in [t, T]$, $s' > s$, we have
\begin{align}\label{xi}
\hE_t[|X^{t,\xi}_s - X^{t,\xi'}_{s'}|^p]
 \leq C\big(|\xi-\xi'|^p+|\xi-\xi'|+(1+|\xi|^p)|s-s'|^{p/2}\big),
\end{align}
and for $\gamma \in [0, T-t]$, we get
\begin{align}\label{Xxi}
\hE_t[\sup_{s\in[t, t+\gamma]}|X^{t,\xi}_s - \xi|^p] \leq C(1+|\xi|^p)\gamma^{\frac{p}{2}}.
\end{align}
\end{remark} 

 We now strict ourselves to a slight modification of the CIR process \eqref{GSDE} as:
\begin{equation}\label{GSDE2}
dX^{t, x}_s= (2\beta_1 X^{t, x}_s+\delta_1)\,ds + (2\beta_2 X^{t, x}_s+\delta_2)\,d\langle B\rangle_s+ \sigma(X^{t, x}_s)\,dB_s, \quad s \geq t,
\end{equation}
where $X^{t, x}_t = x \geq 0$ and $\delta_1, \delta_2$ are positive constants.
\begin{lemma}\label{MarPhi}
Let $T > 0$ be a fixed real number and $X$ holds \eqref{GSDE2}. Then 
\[\hE_t[\varphi(X^{x}_{t+s})] = \hE[\varphi(X^{t,y}_{t+s})]_{y=X^x_t}\] 
for $s \in [0, T-t]$ and $\varphi \in \cliprr$.
\end{lemma}
\begin{proof}
Since the proof for the other cases is similar, we only consider the case that $s=T-t$.  Let $(Y^{t,\xi},Z^{t,\xi},K^{t,\xi})$ be the solution of the following $G$-BSDE
\begin{align}\label{BSDE}
Y^{t,\xi}_s=\varphi(X^{t,\xi}_T)-\int_s^T Z^{t,\xi}_r dB_r-(K^{t,\xi}_T-K^{t,\xi}_s), \quad s \in [t, T]
\end{align}
where $\varphi\in C_{b,Lip}(\mathbb{R})$ 
and $Y^{t,\xi} \in S^{2}_G(0, T)$, $Z^{t, \xi} \in H^2_G(0, T)$ and $K$ is a decreasing $G$-martingale with $K^{t,\xi}_0 = 0$, $K^{t, \xi}_T \in L^2_G(\Omega_T)$. Note that $\xi \in L^2_G(\Omega_t)$.
 By Proposition 2.16 in Hu et al. \cite{HJPS2}, we have
\begin{displaymath}
|Y^{t,\xi}_t-Y^{t,\xi'}_t|^2\leq C\hat{\mathbb{E}}_t[|\varphi(X^{t,\xi}_T)-\varphi(X^{t,\xi'}_T)|^2],
\end{displaymath}
where $C$ is a constant only depending on $T$, $\overline{\sigma}$ and $\underline{\sigma}$. Since $\varphi$ is a bounded Lipschitz, so by \eqref{xi}, we have that
\begin{equation}\label{Yxi}
|Y^{t,\xi}_t-Y^{t,\xi'}_t|^2 \leq C L (|\xi-\xi'|+|\xi-\xi'|^2),
\end{equation}
where $L$ is Lipschitz constant. Now we define
\begin{displaymath}
u(t,x):=Y^{t,x}_t.
\end{displaymath}
Note that $u$ is a deterministic function. 
Thus, $u(t,x)=\hat{\mathbb{E}}[\varphi(X^{t,x}_T]$). It follows from the above estimate that $|u(t,x)-u(t,x')|^2\leq C(|x-x'|+|x-x'|^2)$. By a similar analysis as the proof of Theorem 4.4 in Hu et al. \cite{HJPS2} for any $\xi\in L_G^2(\Omega_t)$, we have
\begin{equation}\label{theorem4.4}
u(t,\xi)=Y^{t,\xi}_t.
\end{equation}
We can write
$\hat{\mathbb{E}}_t[\varphi(X^{x}_T)]
=\hat{\mathbb{E}}_t[\varphi(X^{t,X^{x}_t}_T)]$ and by \eqref{BSDE}, we have $\hat{\mathbb{E}}_t[\varphi(X^{t,X^{x}_t}_T)]=Y^{t,X^x_t}_t$
and due to the definition of  $u$
\begin{align*}
\hat{\mathbb{E}}_t[\varphi(X^{x}_T)]& = u(t,X^x_t)
=\hat{\mathbb{E}}[\varphi(X^{t,y}_T)]|_{y=X^x_t}.
\end{align*}
\end{proof}

Conditional $G$-expectation $\hE_{\tau^{+}}$, where $\tau$ is an optional time, was first constructed by Nutz and Van Handel \cite{NH}  for all upper semianalytic functions. Hu, Ji and Liu \cite{HJL} extended the deterministic-time conditional $G$-expectation notion to optional times for a large class of Borel functions and derived some of regularity properties for $\hE_{\tau^{+}}$ which are helpful
in investigating strong Markov property for a $G$-SDE under the global Lipschitz condition imposed on the coefficients. Here, we proceed the approach used in \cite{HJL} for $G$-SDE \eqref{GSDE1}. One can find the necessary preliminaries and definitions in Appendix. 

We now consider a general one-dimensional $G$-SDE as
\begin{align}\label{GSDEMar}
dX^x_t = \mu(X^x_t)\,dt + b(X^x_t)\,d\langle B\rangle_t+\sigma(X^x_t)\,dB_t, \quad t \in [0, T],
\end{align}
with $X^x_0 = x \in \R$ where $\mu(X_\cdot), b(X_\cdot) \in M_G^1(0,T)$ and $\sigma(X_\cdot) \in M_G^2(0,T)$. We assume that $\mu, b, \sigma$ are sufficiently smooth. Let $\varphi \in {C_{b,Lip}(\mathbb{R}^m)}$ and $0 \leq t_1 \leq t_2 \leq ... \leq t_m <T'$ and $t \geq 0$. The solution $(X^x_t)_{t \geq 0}$ is said to satisfy \emph{Markov property} if 
\begin{align*}
\hE_{t}[\varphi(X^x_{t + t_1},..., X^x_{t + t_m})] = \hE[\varphi(X^y_{t_1},..., X^y_{t_m})]|_{y= X^x_{t}},
\end{align*}
and satisfy \emph{strong Markov property} if
\begin{align*}
\hE_{\tau^{+}}[\varphi(X^x_{\tau + t_1},..., X^x_{\tau + t_m})] = \hE[\varphi(X^y_{t_1},..., X^y_{t_m})]|_{y= X^x_{\tau}}
\end{align*}
for any optional time $\tau$.
\begin{lemma}\label{Mardeter}
Let $\delta^i$, for $i=1,2$, be constant in the $G$-SDE \eqref{GSDE1}.
Then the solution of the $G$-SDE holds Markov property.
\end{lemma}

\begin{proof}
Analogous to the proof of Lemma 4.1 in \cite{HJL}, the proof is carried out by Lemma \ref{MarPhi}.
\end{proof}

We turn to strong Markov property. Thanks to Lemma \ref{pbound}
\[ \mathrm{E}_P[|X^x_t -  X^x_s|^4] \leq  \hE[|X^x_t -  X^x_s|^4] \leq C |t - s|^2 \]
for each $P \in \mathcal{P}$. Then by Kolmogorov's moment criterion for tightness, one can show that 
$\varphi(X^x_{\tau + t_1},..., X^x_{\tau + t_m}) \in L^{1, \tau^{+}}_G(\Omega)$. To more details, see Lemma 4.3 in \cite{HJL}.
\begin{theorem}\label{marcir}
The solution of $G$-SDE \eqref{GSDE1}, under the condition in Lemma \ref{Mardeter}, holds strong Markov property.
\end{theorem}
\begin{proof}
Based on the properties of conditional expectation $\hE_{\tau+}[\cdot]$, Lemma \ref{pbound} and Lemma \ref{Mardeter}, proceeding similarly to the proof of Theorem 4.2 in \cite{HJL}, we obtain the desired result.
\end{proof}

\section{Cox-Ingersoll-Ross model under volatility uncertainty}
In this section, under Knightian uncertainty, 
we focus on the Cox-Ingersoll-Ross model, as a special case of the $G$-SDE \eqref{GSDE}, to describe the instantaneous short rate model as
\begin{align}\label{GCIRmain}
dX_t = (\delta_1-\beta_1 X_t)\,dt + (\delta_2 - \beta_2 X_t)\,d\langle B\rangle_t+ \sigma \sqrt{X_t}\,dB_t, \quad t \geq 0,
\end{align}
where $\delta_1, \delta_2, \beta_1, \beta_2, \sigma$ are positive constants. According to Theorems \ref{main}, \ref{marcir}, $G$-SDE \eqref{GCIRmain} has a unique continuous quasi-surely non-negative solution satisfying strong Markov property.

In the sequel, we aim to inspect the sublinear expectation of the functionals defined on the paths of CIR-type \eqref{GCIRmain}. More specifically, we want to calculate $\hE[\varphi(X_{t'})]$ 
for $t' \in [0, T]$ and some appropriate $\varphi$ such that $\varphi(X_{t'})\in L_G^1(\Omega_{t'})$, where $X$ 
satisfying the dynamics \eqref{GCIRmain}.  This problem have an extensive applications in financial/actuarial markets; for instance payoff of the various products associated with the underlying process. Furthermore, calculating the moments of the underlying process is required to do the necessary analysis in some scientific fields such as statistics and data science.
In order to solve this problem, we need to establish the relation between $\hE[\varphi(X_{t'})]$ and some related PDE. Let $t' \in (0, T]$ be a fixed number. Consider the following forward-backward $G$-SDE:
\smallskip
\begin{align}\label{GBSDE}
	\left\{\begin{array}{lll}
		X^{t,x}_{s} = x + \int_{t}^{s}(\delta_1-\beta_1 X^{t,x}_u)\,du
		+\int_{t}^{s} (\delta_2 - \beta_2 X^{t,x}_u)\, d\langle B\rangle_u + \int_{t}^{s} \sigma \sqrt{X^{t,x}_u}\,dB_u,\\
		Y^{t,x}_s=\varphi(X^{t,x}_{t'})-\int_s^{t'} Z^{t,x}_r dB_r-(K^{t,x}_{t'}-K^{t,x}_s), \qquad 0 \leq t < s \leq t', 
	\end{array}\right.
\end{align}
where $x \geq 0$ and $\varphi \in \cliprr$. 
We denote by $(Y^{t,x},Z^{t,x},K^{t,x})$ the unique solution of the aforementioned $G$-FBSDE \eqref{GBSDE}. Set
\[u^{\varphi}(t,x) := Y^{t,x}_t, \quad (t,x) \in [0, T] \times \R.\] 
Similar to the argument presented in \cite{P} or the proof of Lemma \ref{MarPhi}, we can show that 
\begin{align}\label{Yxxi}
	u^{\varphi}(t,\xi) = Y^{t,\xi}_t
\end{align}
for each $\xi \in L^2_G(\Omega_t)$. Note that $u^{\varphi}$ is a deterministic function as
\begin{align}\label{uvarphi}
	u^{\varphi}(t, x) = \hE[\varphi(X_{t'}^{t,x})].
\end{align}
and using Lemma \ref{MarPhi}
\begin{align}\label{ugamma}
	u^{\varphi}(t,x) = \hE[u^{\varphi}(t+\gamma, X_{t'}^{t+\gamma,X_{t+\gamma}^{t,x}})]
\end{align}
for $\gamma \in [0, t'-t)$. Now, we generalize the \emph{Feynman-Kac} formula to the non-Lipschitz diffusion case. Since $u^\varphi$ may not be a smooth function, we first recall some basic notations and results about viscosity solutions.
Consider the following one-dimensional nonlinear partial differential equation
\begin{align}\label{vis_app}
	\left\{ \begin{array}{lll}
		u_t - F(t, x, u, {u_x}, u_{xx}) =0,& \quad (t,x) \in [0, T] \times \mathbb{R},\\
		u(T, x) = \varphi(x),& \quad x \in \mathbb{R}
	\end{array}\right.
\end{align}
for $\varphi \in C^{0}(\mathbb{R})$ and $F : [0, T] \times \mathbb{R}^4 \rightarrow \mathbb{R}$ a continuous nonlinear function such that
\begin{align}\label{degenerate}
	F(t, x, p, q, a) \geq F(t, x, p, q, b), \quad \textup{if}  \quad a \leq b
\end{align}
for all $t \in [0, T]$, $x, p, q, a, b \in\mathbb{R}$. The relation \eqref{degenerate} indicates that $F$ is $\emph{degenerate elliptic}$.
In general, to solve the fully nonlinear degenerate PDE \eqref{vis_app}, one needs the \emph{viscosity solutions} approach \cite{CIL}. However, in the non-degenerate case, viscosity solution becomes a $C^{1+\frac{\alpha}{2}, 2+\alpha}$ solution for some $\alpha >0$ on $[0,T-\kappa]\times\mathbb{R}$ for some $\kappa >0$, see \cite{CC, P, WL}.

\begin{definition}
Let $u\in C((0,T)\times\mathbb{R})$ and $(t,x)\in(0,T)\times \mathbb{R}$. We denote by $\mathcal{P}^{2,+} u(t,x)$ $($the "parabolic superjet" of $u$ at $(t,x)$$)$ the set of triples $(p,q,X)\in\mathbb{R}^3$ satisfying
\begin{align*}
	u(t+h, x+y)\leq& u(t,x)+ph+qy+\frac{1}{2}Xy^2 + o(|h|+y^2),
\end{align*}
where $(t+h , x+y) \in (0,T)\times \mathbb{R}$ and $h, y \rightarrow 0$.
Similarly, we define $\mathcal{P}^{2,-} u(t,x)$ $($the "parabolic subjet" of $u$ at $(t,x) )$ by $\mathcal{P}^{2,-} u(t,x):=-\mathcal{P}^{2,+}(- u)(t,x)$.
\end{definition}

Now we define the viscosity solution.
\begin{definition}\label{viscositysolution}
\cite{CIL, KN}. Let $\Lambda := [0, T]\times\mathbb{R}$. The function $u\in C(\Lambda)$ is called a viscosity subsolution of \eqref{vis_app} at $(t,x)\in\Lambda$ if
\begin{displaymath}
\inf_{(p,q,X)\in\mathcal{P}^{2,+}u(t,x)} p - F(t, x, u(t, x), q, X) \geq 0.
\end{displaymath}
The function $u\in C(\Lambda)$ is called a viscosity supersolution of \eqref{vis_app} at $(t,x)\in\Lambda$ if 
\begin{displaymath}
	\sup_{(p,q,X)\in\mathcal{P}^{2,-}u(t,x)} p - F(t, x, u(t, x), q, X) \leq 0.
\end{displaymath}
	If $\mathcal{P}^{2,+}u(t,x) = \emptyset$, then
	\begin{displaymath}
		\inf_{(p,q,X)\in\mathcal{P}^{2,+}u(t,x)} p - F(t, x, u(t, x), q, X) = \infty.
	\end{displaymath}
	A similar result is obtained for the case of $\mathcal{P}^{2,-}u(t,x)= \emptyset$.
\end{definition}

By remembering definition of $u^\varphi$ in \eqref{uvarphi}, we present a new result:
\begin{theorem}\label{fey}
	The function $u^\varphi$ is the unique viscosity solution of the following nonlinear parabolic PDE:
	\smallskip
	\begin{align}\label{PDEnolin}
		\left\{ \begin{array}{lll}
			u_t(t,x) + (\delta_1-\beta_1 x)u_x(t,x) + 2\,G((\delta_2-\beta_2 x)u_x(t,x) + \frac{\sigma^2 x}{2}  u_{xx}(t,x))=0,\\
			(t,x) \in [0, t') \times \R^{+},\\
			u(t', x)=\varphi(x),\quad x \in \R^{+}.
		\end{array}\right.
	\end{align}
	where $G$ has already been defined in Section 2. 
\end{theorem}
\begin{proof}
The proof is analogous to the one of the nonlinear Feynman-Kac theorem to the Lipschitz $G$-SDE by Peng \cite{P} or Hu et al. \cite{HJPS2}.
\end{proof} 

Utilizing the Theorem \ref{fey}, analogous to Yang and Zhao in \cite{JYWZ}, letting $\varphi(x) := I_{\{x \leq a\}}$, $a \in \R$, and setting $F_{X^{t,x}_{t'}}(a) := \hE[I_{\{X^{t,x}_{t'}\leq a\}}]$ and simulating $u^{\varphi}(t,x)$ by efficient numerical methods, see e.g., \cite{GPE, MCLR, SPE}, we obtain the distribution of $X^{t,x}_{t'}$ by \eqref{uvarphi}. 
Similarly, by setting $\varphi(x) = x$, $\varphi(x) = -x$, $\varphi(x) = x^2$ and $\varphi(x)=-x^2$ and computing $u^{\varphi}$ from \eqref{PDEnolin}, mean-uncertainty and variance-uncertainty of $X^{t,x}_{t'}$ can be obtained, respectively.

\begin{proposition}\label{moments of X}
For the case that $\delta_2 = \beta_2 = 0$ and $\beta_1\neq 0$, we have
\begin{equation}\label{nomeanuncertainty}
\hat{\mathbb{E}}[X^{t,x}_{t'}]=-\hat{\mathbb{E}}[-X^{t,x}_{t'}]=e^{-\beta_1 (t'-t)}(x-\frac{\delta_1}{\beta_1})+\frac{\delta_1}{\beta_1},
\end{equation}
and 
\begin{align}\label{x^2}
\hat{\mathbb{E}}[|X^{t,x}_{t'}|^2]=u(t,x;t',\overline{\sigma}), \ -\hat{\mathbb{E}}[-|X^{t,x}_{t'}|^2]=u(t,x;t',\underline{\sigma}),
\end{align}
where 
\begin{equation}\begin{split}\label{usigma}
u(t,x;t',a)=&\frac{\delta^2_1}{\beta^2_1}+\frac{\delta_1\sigma^2 a^2}{2\beta^2_1}+e^{-\beta_1(t'-t)}\frac{\sigma^2a^2+2\delta_1}{\beta_1}(x-\frac{\delta_1}{\beta_1})\\
&+e^{-2\beta_1(t'-t)}\{(x-\frac{\delta_1}{\beta_1})^2+\frac{\sigma^2a^2}{2\beta_1}(\frac{\delta_1}{\beta_1}-2x)\}.
\end{split}\end{equation}

For the case that $\delta_1 = \beta_1 = 0$ and $\beta_2\neq 0$, we have
\begin{align}\label{meanuncertainty}
\hE[X^{t,x}_{t'}]=u(t,x), \ \hE[-X^{t,x}_{t'}]=\tilde{u}(t,x),
\end{align}
where
\begin{align}\label{uu}
u(t,x)=\begin{cases}
(x-\frac{\delta_2}{\beta_2})e^{\overline{\sigma}^2\beta_2(t-t')}+\frac{\delta_2}{\beta_2}, & 0\leq x\leq\frac{\delta_2}{\beta_2},\\[4mm]
(x -\frac{\delta_2}{\beta_2})e^{\underline{\sigma}^2\beta_2(t-t')}+\frac{\delta_2}{\beta_2}, & x\geq\frac{\delta_2}{\beta_2}, 
\end{cases}
\end{align}
and 
\begin{align}\label{uu'}
\tilde{u}(t,x)=\begin{cases}
(\frac{\delta_2}{\beta_2}-x)e^{\underline{\sigma}^2\beta_2(t-t')}-\frac{\delta_2}{\beta_2}, & 0\leq x\leq\frac{\delta_2}{\beta_2},\\[4mm]
(\frac{\delta_2}{\beta_2} - x)e^{\overline{\sigma}^2\beta_2(t-t')}-\frac{\delta_2}{\beta_2}, & x\geq\frac{\delta_2}{\beta_2},
\end{cases}
\end{align}
for $(t,x) \in [0,t'] \times \mathbb{R}^+$.
\end{proposition}

\begin{proof}
First, we consider the case $\delta_2 = \beta_2 = 0$, $\beta_1\neq 0$. 
Applying $G$-It\^{o}'s formula to $e^{\beta_1 s}X^{t,x}_s$, it is easy to get that 
\begin{displaymath}
e^{\beta_1 t'}X^{t,x}_{t'}=e^{\beta_1 t}x+\int_t^{t'} \delta_1 e^{\beta_1 s}ds+\int_t^{t'} \sigma e^{\beta_1s}\sqrt{X^{t,x}_s}dB_s.
\end{displaymath}
Noting that $\int_t^\cdot \sigma e^{\beta_1s}\sqrt{X^{t,x}_s}dB_s$ is a symmetric $G$-martingale and taking expectations on both sides of the above equation, we finally obtain that 
\begin{align*}
\hat{\mathbb{E}}[X^{t,x}_{t'}]=-\hat{\mathbb{E}}[-X^{t,x}_{t'}]=e^{-\beta_1 (t'-t)}(x-\frac{\delta_1}{\beta_1})+\frac{\delta_1}{\beta_1}.
\end{align*}
By Equation \eqref{usigma}, 
it is easy to check that $xu_{xx}(t,x;t',\overline{\sigma}) \geq 0$ and the function $u(\cdot,\cdot;t',\overline{\sigma})$ satisfies 
\begin{align*}
u_t(t,x) + (\delta_1 -\beta_1x)u_x(t,x)
+ \frac{\overline{\sigma}^2\sigma^2}{2} x\,u_{xx}(t,x)&=0, \quad u(t',x) = x^2.
\end{align*}
Therefore, it is the solution of PDE \eqref{PDEnolin} with $\delta_2 = \beta_2 = 0$ and $\varphi(x)= x^2$. By Theorem \ref{fey}, we have
\begin{align*}
\hat{\mathbb{E}}[|X^{t,x}_{t'}|^2]=u(t,x;t',\overline{\sigma}).
\end{align*}
In a similar way, we can obtain that 
\begin{align*} 
-\hat{\mathbb{E}}[-|X^{t,x}_{t'}|^2]=u(t,x;t',\underline{\sigma}).
\end{align*}

Now we compute $\hE[X^{t,x}_{t'}]$ and $\hE[-X^{t,x}_{t'}]$ for the case that $\delta_1=\beta_1=0$, $\beta_2 \neq 0$. By Theorem \ref{fey}, it suffices to check that $u$ and $\tilde{u}$ defined by \eqref{uu} and \eqref{uu'} are the viscosity soutions of the following PDEs
\begin{align}\label{uvar}
u_t(t,x) + 2\,G((\delta_2-\beta_2 x)u_x(t,x) + \frac{\sigma^2 x}{2}  u_{xx}(t,x))& = 0, \quad u(t', x) = x.
\end{align}
and
\begin{align}\label{uvar'}
	\tilde{u}_t(t,x) + 2\,G((\delta_2-\beta_2 x)\tilde{u}_x(t,x) + \frac{\sigma^2 x}{2}  \tilde{u}_{xx}(t,x))& = 0, \quad \tilde{u}(t', x) = -x.
\end{align}
respectively. We only prove the case for $u$. It is easy to check that $u$ is the classical solution of the PDE for $\{(t,x) \in  [0, t'] \times [0, \frac{\delta_2}{\beta_2}) \}$
and $\{(t,x) \in  [0, t'] \times (\frac{\delta_2}{\beta_2}, \infty)\}$. For the area $\{(t,x)\in[0,t']\times\{\frac{\delta_2}{\beta_2}\}\}$ under which $u$ is not differentiable, we need to compute $\mathcal{P}^{2,-}u(t,\delta_2/\beta_2) $ and $\mathcal{P}^{2,+}u(t,\delta_2/\beta_2) $. We claim that 
\begin{align}\label{P-}
\mathcal{P}^{2,-}u(t,\delta_2/\beta_2)=\{(p,q,X): (p,q,X)\in (-\infty,0] \times [e^{\overline{\sigma}^2\beta_2(t-t')}, e^{\underline{\sigma}^2\beta_2(t-t')}]\times (-\infty, 0]\}
\end{align}
and 
\begin{align}\label{P+}
\mathcal{P}^{2,+}u(t,\delta_2/\beta_2) = \emptyset.
\end{align}
Then, by Definition \ref{viscositysolution}, $u$ defined by \eqref{uu} is indeed the viscosity solution to PDE \eqref{uvar}.

In the following, we prove \eqref{P-} and \eqref{P+} hold. Note that 
\begin{align*}
\mathcal{P}^{2,-}u(t,\frac{\delta_2}{\beta_2}) & = \{(p, q, X) \,|\, u(t+h, \frac{\delta_2}{\beta_2}+ y) \geq  u(t,\frac{\delta_2}{\beta_2}) + ph + qy + \frac{1}{2}Xy^2 + o(|h|+y^2) \},
\end{align*}
where $h, y \rightarrow 0$. Hence, for any $(p,q,X)\in \mathcal{P}^{2,-}u(t,\delta_2/\beta_2)$, we have
\begin{align*}
p+2G((\delta_2-\beta_2 \frac{\delta_2}{\beta_2})q+\frac{1}{2}\sigma^2\frac{\delta_2}{\beta_2} X)\leq 0.
\end{align*} 
Also, we have 
\begin{align}\label{uh}
\frac{u(t + h, \delta_2/\beta_2 + y) -  u(t, \delta_2/\beta_2)}{h} - p \geq (qy + \frac{1}{2}Xy^2 + o(|h| + y^2))/h.
\end{align}
It is easy to check that for $y > 0$
\begin{align*}
\lim_{h\rightarrow 0}\lim_{y\rightarrow 0}\frac{u(t + h, \delta_2/\beta_2 + y) -  u(t, \delta_2/\beta_2)}{h}=\lim_{h\rightarrow 0}\lim_{y\rightarrow 0}\frac{y e^{\underline{\sigma}^2\beta_2(t+h-t')}}{h}
&=0
\end{align*}
and for $y<0$
\begin{align*}
\lim_{h\rightarrow 0}\lim_{y\rightarrow 0}\frac{u(t + h, \delta_2/\beta_2 + y) -  u(t, \delta_2/\beta_2)}{h}=\lim_{h\rightarrow 0}\lim_{y\rightarrow 0}\frac{y e^{\overline{\sigma}^2\beta_2(t+h-t')}}{h}
&=0.
\end{align*}
So by \eqref{uh}, tending $y\rightarrow 0$ and then $h\rightarrow 0$ results in $p \leq 0$. 

Similarly, for $y>0$, we have
\begin{align*}
\frac{u(t + h, \delta_2/\beta_2 + y) -  u(t, \delta_2/\beta_2)}{y} - q \geq (ph + \frac{1}{2}Xy^2 + o(|h| + y^2))/y.
\end{align*}
By tending $h, y \rightarrow 0$, we get $q \leq  e^{\underline{\sigma}^2\beta_2(t-t')}$. For $y < 0$, we get
\begin{align*}
\frac{u(t + h, \delta_2/\beta_2 + y) -  u(t, \delta_2/\beta_2)}{y} - q \leq (ph + \frac{1}{2}Xy^2 + o(|h| + y^2))/y.
\end{align*}
By tending $h, y \rightarrow 0$, we get $q \geq  e^{\overline{\sigma}^2\beta_2(t-t')}$. 

Furthermore, for $y<0$, we can write
\begin{align*}
y e^{\overline{\sigma}^2\beta_2(t+h-t')}\geq  ph + qy + \frac{1}{2}Xy^2 + o(|h|+y^2).
\end{align*}
So by tending $h \rightarrow 0$, we have
\begin{align*}
y e^{\overline{\sigma}^2\beta_2(t-t')} -qy &\geq \frac{1}{2}Xy^2 + o(y^2).
\end{align*}
Similarly for $y>0$ and by tending $h \rightarrow 0$, we get
\begin{align*}
y e^{\underline{\sigma}^2\beta_2(t-t')} - qy \geq  \frac{1}{2}Xy^2 + o(y^2).
\end{align*}
These result in $X \in (-\infty, 0]$. Finally, we have 
\begin{align*}
\mathcal{P}^{2,-}u(t,\delta_2/\beta_2)=\{(p,q,X): (p,q,X)\in (-\infty,0] \times [e^{\overline{\sigma}^2\beta_2(t-t')}, e^{\underline{\sigma}^2\beta_2(t-t')}]\times (-\infty, 0]\}.
\end{align*}



For $\mathcal{P}^{2,+}u(t,\delta_2/\beta_2)$, we see that 
\begin{align*}
\frac{u(t + h, \delta_2/\beta_2 + y) -  u(t, \delta_2/\beta_2)}{y} - q \leq (ph + \frac{1}{2}Xy^2 + o(|h| + y^2))/y
\end{align*}
for $y>0$. By tending $h, y \rightarrow 0$, we get $q \geq  e^{\underline{\sigma}^2\beta_2(t-t')}$. For $y < 0$, we get 
\begin{align*}
\frac{u(t + h, \delta_2/\beta_2 + y) -  u(t, \delta_2/\beta_2)}{y} - q \geq (ph + \frac{1}{2}Xy^2 + o(|h| + y^2))/y.
\end{align*}
By tending $h, y \rightarrow 0$, we get $q \leq  e^{\overline{\sigma}^2\beta_2(t-t')}.$ This results in that $\mathcal{P}^{2,+}u(t,\delta_2/\beta_2) = \emptyset$.
\end{proof}
\begin{remark}
If $\delta_2 = \beta_2 = 0$, i.e. there is no $d\langle B\rangle$-term, the moments of the first and second order of the CIR process under $G$-Brownian motion in \eqref{nomeanuncertainty} and \eqref{x^2} are exactly equal to those of the classical CIR process under Brownian motion by a slight difference under which $\sigma$ is replaced by $\sigma \overline{\sigma}$.  However, $\hat{\mathbb{E}}[|X^{t,x}_{t'}|^2]$ and $-\hat{\mathbb{E}}[-|X^{t,x}_{t'}|^2]$ are not equal under $G$-expectation while they are clearly equal under classical expectation. If $\beta_2\neq 0$, i.e., the quadratic variation term appears, there is mean uncertainty for the $G$-CIR process.
\end{remark}

\section*{Conclusion}
This paper extends the famous Cox-Ingersoll-Ross process to the case of ambiguity which is modeled by $G$-Brownian motion. 
The existence and uniqueness of the solution of the aforementioned CIR process are accomplished. Then, by estimating the process and also the BSDE approach, the strong Markov property of the process is investigated. In the final section, based on the extended nonlinear Feynman-Kac theorem to the case of underlying CIR process, distributions of the CIR process, i.e., the expectation of functionals of the process, are provided. Then, we estimate some moments of the process for some special parameters in the equation. This paper is the first contribution to the CIR process driven by $G$-Brownain motion and we believe that it may serve as an inspiration for future studies from both mathematical and financial aspects.

\section*{Acknowledgments}
 Li's research was supported by the Natural Science Foundation of Shandong Province (No. ZR2022QA022) and the Qilu Young Scholars Program of Shandong University.
 
 \renewcommand\thesection{Conflict of Interest}
\section{ }
 The authors declared that they have no conflict of interest.

\appendix
\renewcommand\thesection{Appendix}
\section{ }
\label{appA}
\renewcommand\thesection{A}
\subsection{Conditional $G$-expectation $\hE_{\tau^{+}}$}
In this section, we review some preliminaries on Conditional $G$-expectation $\hE_{\tau^{+}}$ for any given optional time $\tau$.
To do this, first, set
\begin{itemize}
\item $\mathcal{F}_t:= \sigma(\{B_s: s\leq t\})$
\item $L^0(\mathcal{F}_t): = \{X: X\,\, \textup{is}\,\, \mathcal{F}_t$-$\textup{measurable}\}$
\item $\mathcal{F}:= {\bigvee_{t \geq 0}}\mathcal{F}_t$
\item $L^0(\mathcal{F}): = \{X: X\,\, \textup{is}\,\, \mathcal{F}$-$\textup{measurable}\}$
\end{itemize}
\begin{definition}
The non-negative random variable $\tau$ $(\tau: \Omega \rightarrow [0, \infty))$ is called a stopping time if $\{\tau \leq t\} \in \mathcal{F}_t$ for each $t\geq 0$ and an optional time if $\{\tau < t\} \in \mathcal{F}_t$ for each $t \geq 0$.
\end{definition}

For each optional time $\tau$, we set the following $\sigma$-field
\begin{equation*}
\mathcal{F}_{\tau^+} :=\{A \in \mathcal{F} : A \cap \{\tau < t\} \in \mathcal{F}_t, \forall\, t \geq 0\}= \{A \in \mathcal{F} : A \cap \{\tau \leq t\} \in \mathcal{F}_{t^+}, \forall\, t \geq 0\},
\end{equation*}
where $\mathcal{F}_{t^+} = \cap _{s>t} \mathcal{F}_s$. In the case of stopping time, we also set
\begin{equation*}
\mathcal{F}_{\tau} :=\{A \in \mathcal{F} : A \cap \{\tau \leq t\} \in \mathcal{F}_t, \forall\, t \geq 0\}.
\end{equation*}
For optional time $\tau$ and $p \geq 1$, we set
\begin{displaymath}
L^{0,p,\tau^+}_G(\Omega) = \{X = \sum_{i=1}^{n}
\xi_i I_{A_i} : n \in \mathbb{N}, \{A_i\}_{i=1}^n \textup{is an}\, \mathcal{F}_{{\tau}^+}\textup{-partition of}\,\, \Omega,\, \xi_i \in L_G^p(\Omega), i = 1,... ,n\}.
\end{displaymath}
We denote by $L^{p,\tau^+}_G(\Omega)$ the completion of $L^{0,p,\tau^+}_G(\Omega)$ under the norm $\|.\|_p$.
The authors \cite{HJL} constructed conditional $G$-expectation $\hE_{\tau^+}$ on $Lip(\Omega)$, $L^1_G(\Omega)$ and $L^{1,\tau^+}_G(\Omega)$  in three stages.
\begin{definition}
For $X \in Lip(\Omega)$, we define $\hE_{\tau^+} : Lip(X) \rightarrow L^{1,\tau^+}_G(\Omega) \cap L^0(\mathcal{F}_{\tau^+})$ as
\begin{itemize}
\item Let $\tau$ be a simple discrete stopping time $\tau$ taking values in $\{t_i : i \geq 1\}$, \\
\[\hE_{\tau^+}[X] := \sum_{i=1}^{\infty} \hE_{t_i}[X] I_{\{\tau = t_i\}},\]
where $t_i \uparrow \infty$ as $i \rightarrow \infty$. We employ the convection that $t_{n+i}:= t_n + i$, $i \geq 1$, if $\tau$ is a discrete stopping time taking finitely many values $\{t_i: i \leq n\}$ with $t_i \leq t_{i+1}$.
\item For a general optional time $\tau$, the conditional $G$-expectation is defined by 
\[\hE_{\tau^+}[X] := \lim_{n \rightarrow \infty} \hE_{{\tau_n}^+} [X] \quad in \quad \mathbb{L}^1,\]
where $\tau_n$ be a sequence of simple discrete stopping times such that $\tau_n \rightarrow \tau$ uniformly.
\end{itemize}
\end{definition}
\begin{definition}
For $X \in L^1_G(\Omega)$, we define $\hE_{\tau^+} : L^1_G(\Omega) \rightarrow L^{1,\tau^+}_G(\Omega) \cap L^0(\mathcal{F}_{\tau^+})$ as
\[\hE_{\tau^+}[X] := \lim_{n \rightarrow \infty} \hE_{{\tau}^+} [X_n] \quad in \quad \mathbb{L}^1,\]
where $\{X_n\}_{n=1}^{\infty} \subset Lip(\Omega)$ such that $X_n \rightarrow X$ in $\mathbb{L}^1$.
\end{definition}
\begin{definition}
We define $\hE_{\tau^+} : L^{1, \tau^{+}}_G(\Omega) \rightarrow L^{1,\tau^+}_G(\Omega) \cap L^0(\mathcal{F}_{\tau^+})$ in two steps. Let $X = \sum_{i=1}^{n} \xi_i I_{A_i} \in L^{0,1,\tau^{+}}_G(\Omega)$ where $\xi_i \in L^1_G(\Omega)$
and $\{A_i\}_{i=1}^{n}$ is an $\mathcal{F}_{\tau^+}$-partition of $\Omega$. We define
\[\hE_{\tau^+}[X] := \sum_{i=1}^{\infty} \hE_{\tau^+}[\xi_i] I_{A_i},\]
Let $X \in L^{1, \tau^{+}}_G(\Omega)$. Then we define
\[\hE_{\tau^+}[X] := \lim_{n \rightarrow \infty} \hE_{{\tau}^+} [X_n] \quad in \quad \mathbb{L}^1,\]
where $\{X_n\}_{n=1}^{\infty} \subset L^{0,1,\tau^{+}}_G(\Omega)$ such that $X_n \rightarrow X$ in $\mathbb{L}^1$.
\end{definition}
\begin{proposition}\label{conditionalexpectation}
The conditional $G$-expectation $\hE_{\tau^+}$ defined through three definitions above is well-defined and for $X, Y \in L^{1, \tau^+}_G$
\begin{itemize}
\item $\hE_{\tau^+}[X] \leq \hE_{\tau^+}[Y]$, for $X \leq Y$;
\item $\hE_{\tau^+}[X+Y] \leq \hE_{\tau^+}[X] + \hE_{\tau^+}[Y]$;
\item $\hE[\hE_{\tau^+}[X]] = \hE[X]$.
\end{itemize}
\end{proposition}

Note that the Proposition above can be presented for $Lip(\Omega)$ and $L^{0, 1, \tau^+}_G$separately. For more properties of $\hE_{\tau^+}$, the reader can refer to \cite{HJL}. 

\subsection{BSDE driven by  $G$-Brownian motion}
In the section, we present a short discussion on the literature of backward stochastic differential equation driven by $G$-Brownian motion $(B_t)_{t \geq 0}$ whose solution  consists of a triple of processes $(Y, Z, K)$ as:
\begin{align}\label{backward1}
Y_s&= \eta + \int_{s}^{T}f(v,, Y_v, Z_v)\,dv + \int_{s}^{T} g(v, Y_v, Z_v)\, d\langle B\rangle_v-\int_{s}^{T} Z_v\,dB_v -(K_T - K_s),
\end{align}
for $0 \leq s \leq T$ where $K$ is a decreasing $G$-martingale and
$f, g :[0,T] \times \Omega_T \times \R^2 \rightarrow \R$
satisfies the following properties:  
\begin{itemize}
	\item[(Hi)] 
		For all $\omega \in \Omega$, $t \in [0, T]$ and $y', y, z, z' \in \mathbb{R}$
		\begin{align}\label{Lip}
		|f(t, \omega, y, z) - f(t, \omega, y', z')| + |g(t, \omega, y, z) - g(t, \omega, y', z')| \leq C\big[|y - y'| +|z - z'| \big]
	\end{align}
 where $C$ is a positive constant;
	\item[(Hii)] for any $y,z\in\mathbb{R}$,  $f(.,.,y,z)$ and $g(.,.,y,z) \in  M_G^\beta(0, T)$ with some $\beta>1$.
\end{itemize}
For simplicity, we denote by $\mathfrak{S}^{\alpha}(0,T)$ the collection of the processes $(Y, Z, K)$ such that $Y \in S^{\alpha}_G(0,T)$ and $Z \in H^{\alpha}_G(0,T)$ and $K$ is a decreasing $G$-martingale with $K_0 = 0$ and $K_T \in L^{\alpha}_G(\Omega_T)$. 

\begin{definition}\label{defBSDE}
 Let $\eta \in L_G^{\beta}(\Omega_T)$ for some $\beta > 1$. A triple $(Y, Z, K)$ is called a solution of the equation \eqref{backward1} if for some $1 < \alpha \leq \beta$, we have that $(Y, Z, K) \in \mathfrak{S}^{\alpha}(0, T)$ and $(Y,Z,K)$ 
 satisfies equation \eqref{backward1}.
\end{definition}
\begin{theorem}\label{exiBSDE}
Assume that $\eta \in L_G^{\beta}(\Omega_T)$ for some $\beta > 1$ and conditions (Hi), (Hii) are satisfied. Then equation \eqref{backward1} has a unique solution $(Y, Z, K) \in \mathfrak{S}^{\alpha}(0, T)$ for any $1 < \alpha < \beta$.
\end{theorem}

The reader can refer to \cite{HJPS1, HJPS2} for more detailed discussions.





\begin{thebibliography}{00}



\bibitem{ABMO}
\textsc{Akhtari, B., Biagini, F., Mazzon, A., Oberpriller. K.}\ (2020).
Generalized Feynman-Kac Formula under volatility uncertainty.
\emph{arxiv.org/abs/2012.08163.}


\bibitem{BL}
\textsc{Bai, X., Lin, Y.} \ (2014).
On the existence and uniqueness of solutions to stochastic differential equations driven by $G$-Brownian motion with integral-Lipschitz coefficients. 
\emph{Acta Math. Appl. Sin. Engl. Ser.} \textbf{30}, 589-610.


\bibitem{GPE}
\textsc{Barles, G. and Souganidis, P.E.}\ (1991). 
Convergence of approximation schemes for fully nonlinear second order equations.
\emph{Asymptotic analysis} 
\textbf{4(3)}, 271-283.



\bibitem{SBKN}
\textsc{Biagini, S., Bouchard, B., Kardaras, C., Nutz, M.} \ (2017).
Robust fundamental theorem for continuous processes.
\emph{Mathematical Finance}
\textbf{27(4)},
 963-987.



\bibitem{MCLR}
\textsc{Briani, M., Chioma, C.L. and Natalini, R.}\ (2004). 
Convergence of numerical schemes for viscosity solutions to integro-differential degenerate parabolic problems arising in financial theory.
\emph{Numerische Mathematik} 
\textbf{98(4)}, 607-646.


\bibitem{BA} 
\textsc{Brigo, D., Alfonsi, A.}\ (2005).
Credit default swap calibration and derivatives pricing with the SSRD stochastic intensity model.
\emph{Finance and stochastics}
\textbf{9(1)}, 29–42. 


\bibitem{BM} \textsc{Brigo, D., Mercurio, F.}\ (2001).
A deterministic-shift extension of analytically–tractable and time–homogeneous short–rate models.
\emph{Finance and Stochastics} \textbf{5(3)}, 369-387.


\bibitem{BD} \textsc{Brown, Stephen J., Dybvig, Philip H.}\ (1986).
The empirical implications of the Cox, Ingersoll, Ross theory of the term structure of interest rates.
\emph{The Journal of Finance, Wiley Online Library} \textbf{41(3)}, 617-630.


\bibitem{CC}
\textsc{Cabre, X., Caffarelli L. A.}\ (1997).
Fully nonlinear elliptic partial differential equations.
\emph{American Math. Society.}



\bibitem{CIRS1}
\textsc{Cox, J. C., Ingersoll, J. E., Ross, S. A.}\ (1985). 
An intertemporal general equilibrium model of asset prices.
\emph{Econometrica: Journal of the Econometric Society} 363-384. 


\bibitem{CIRS}
\textsc{Cox, J. C., Ingersoll, J. E., Ross, S. A.}\ (1985). 
A theory of the term structure of interest rates.
\emph{Econometrica} 
\textbf{53(2)}, 385-407. 


\bibitem{CIL} \textsc{Crandall, M., Ishii, H., Lions P. L.}\ (1992).
User's guide to the viscosity solutions of second order partial differential equations. 
\emph{Bull. Amer. Math. Soc.}
\textbf{27}, 1-67.


\bibitem{DD} \textsc{Deelstra, G., Delbaen, F.}\ 
Existence of solutions of stochastic differential equations related to the Bessel process. Available at \textrm{ https://people.math.ethz.ch/~delbaen/ftp/preprints/existentie$\_$rev.pdf}

\bibitem {DHP11} \textsc{Denis, L., Hu, M., Peng, S.}\ (2011).
 Function spaces and capacity related to a sublinear expectation: application to $G$-Brownian motion paths. 
 \emph{Potential Anal.} 
 \textbf{34}, 139-161.

\bibitem{EJ13}
{Epstein, L., Ji, S.} \ (2013).
Ambiguous volatility and asset pricing in continuous time.
\emph{The Review of Financial Studies} 
\textbf{26(7)}, 1740-1786.



\bibitem{EJ14}
{Epstein, L., Ji, S.} \ (2014).
Ambiguous volatility, possibility and utility in continuous time.
\emph{Journal of Mathematical Economics} 
\textbf{50}, 269-282. 

\bibitem{FNS}
\textsc{Fadina, T., Neufeld, A.,  Schmidt, T.} \ (2019).
Affine processes under parameter uncertainty.
\emph{Probability, uncertainty and quantitative risk}
\textbf{4(1)}, 1-35.


\bibitem{G}
\textsc{Gao, F.}\ (2009).
Pathwise properties and homeomorphic flows for stochastic differential equations driven by G-Brownian.
\emph{Stochastic Processes and their Applications}
\textbf{119}, 3356–3382.


\bibitem{GR} \textsc{Gy\"{o}ngy, I., R\'{a}sonyi, M.}\ (2011). 
A note on Euler approximations for SDEs with H\"{o}lder continuous diffusion coefficients. \emph{Stochastic Prccesses and their Applications} \textbf{121}, 2189-2200.



\bibitem{H1}
\textsc{H\"{o}lzermann, J.}\ (2021). 
The Hull-White model under volatility uncertainty.
\emph{Quantitative Finance} \textbf{21(11)}, 1921-1933.



\bibitem{H2}
\textsc{H\"{o}lzermann, J.}\ (2020). 
Pricing interest rate derivatives under volatility uncertainty.
\emph{arxiv.org/abs/2003.04606.}


\bibitem{H3}
\textsc{H\"{o}lzermann, J.}\ (2022). 
Term structure modeling under volatility uncertainty.
\emph{Mathematics and Financial Economics}
\textbf{16}, 317-343.



\bibitem{HJL}
\textsc{Hu, M., Ji, X., Liu, G.}\, (2021).
On the strong Markov property for stochastic differential equations driven by G-Brownian motion.
\emph{Stochastic Processes and their Applications}
\textbf{131}, 417-453.



\bibitem{HJPS1}
\textsc{Hu, M., Ji, S., Peng, S., Song, Y.}\, (2014).
Backward stochastic differential equations driven by G-Brownian motion.
\emph{Stochastic Processes and their Applications}
\textbf{124(1)}, 759-784.
 
  
\bibitem{HJPS2}
\textsc{Hu, M., Ji, S., Peng, S., Song, Y.}\, (2014).
Comparison theorem, Feynman-Kac formula and Girsanov transformation for BSDEs driven by G-Brownian motion.
\emph{Stochastic Processes and their Applications}
\textbf{ 124(2)}, 1170-1195.

\bibitem{HWZ} 
\textsc{Hu, M., Wang, F., Zheng, G.}\ (2016)
 Quasi-continuous random variables and processes under the $G$-expectation framework. 
 \emph{Stochastic Process. Appl.} 
 \textbf{126(8)}, 2367-2387.

\bibitem{KS} \textsc{Karatzas, I., Shreve, S.E.}\ (1988). 
\emph{Brownian Motion and Stochastic Calculus}. Springer. New York.



\bibitem{KN}
\textsc{Katzourakis, Nikos.}\ (2014). 
\emph{An Introduction to Viscosity Solutions for Fully Nonlinear PDE with Applications to Calculus of Variations in $L^{\infty}$}. Springer. 



\bibitem{KFH}
\textsc{Knight, F. H.} \ (1921). 
\emph{Risk, uncertainty and profit.}
\textbf{31}, Houghton Mifflin.
 

\bibitem{MY}
\textsc{Maghsoodi, Y.} \ (1996).
Solution of the extended CIR term structure and bond option valuation.
\emph{Math. Finance.}
\textbf{6(1)}, 89-109.



\bibitem{NH}
\textsc{Nutz, M., Van Handel, R.} \ (2013).
Constructing sublinear expectations on path space.
\emph{Stochastic Processes and their Applications.}
\textbf{123(8)}, 3100-3121.


\bibitem{OMB}
\textsc{Orlando, G., Mininni, R. M., Bufalo, M.}\ (2020). 
Forecasting interest rates through Vasicek and CIR models: A partitioning approach.
\emph{Journal of Forecasting}
\textbf{39(4)}, 569-579.


\bibitem{P0}
\textsc{Peng, S.}\ (2007). 
$G$-Brownian Motion and Dynamic Risk
Measure under Volatility Uncertainty.
\emph{arXiv:0711.2834v1.}



\bibitem{P}
\textsc{Peng, S.}\ (2019). 
\emph{Nonlinear Expectations and Stochastic Calculus Under Uncertainty: With Robust CLT and G-Brownian Motion.}
Probability Theory and Stochastic Modelling 95. Springer.



\bibitem{SPE}
\textsc{Souganidis, P. E.}\ (1985). 
Approximation schemes for viscosity solutions of Hamilton-Jacobi equations.
 \emph{Journal of differential equations} 
\textbf{59(1)}, 1-43.
 
 
 
\bibitem{WL}
\textsc{Wang, L.} \ (1992).
On the regularity of fully nonlinear parabolic equations: II.
\emph{Comm. Pure Appl. Math.}
\textbf{45}, 141-178.



\bibitem{Y} \textsc{Yamada, T.}\ (1978). 
Sur une construction des solutions d'\'{e}quations diff\'{e}rentielles stochastiques dans le cas non-Lipschitzien, in: S\'{e}m de Prob. XII, LNM 649. Springer-Verlag. Berlin. Heidelberg, New York. 114-131.



\bibitem{JYWZ}
\textsc{Yang, J., Zhao, W.}\  (2016).
Numerical simulations for G-Brownian motion.
\emph{ Frontiers of Mathematics in China.}
\textbf{11}, 1625-1643.



\end{thebibliography}


\end{document}